\def \version {25 January, 2018}

\documentclass[12pt]{article}
\usepackage{amsmath,amsthm,amssymb}

\def \Rs {\mbox{\sf Rs}}
\newcommand{\R}[1]{{\mbox{\sf R}(#1)}}

\def \Ts {\mbox{\sf Ts}}
\def \ndR {\mbox{\sf R}}
\def \cF {{\cal F}}
\def \cG {{\cal G}}
\def \sst {{\subset}}
\def \ssq {{\subseteq}}
\def \smin {{\setminus}}

\def \nin {\noindent}
\def \bsk {\bigskip}
\def \msk {\medskip}
\def \pf {\nin{\bf Proof.} \ }
\def \qed {\hfill $\Box$}
\def \ov {\overline}

\newtheorem{Theorem}{Theorem}
\newtheorem{lem}[Theorem]{Lemma}
\newtheorem{defi}[Theorem]{Definition}
\newtheorem{crl}[Theorem]{Corollary}
\newtheorem{prp}[Theorem]{Proposition}
\newtheorem{prm}[Theorem]{Problem}
\newtheorem{rmk}[Theorem]{Remark}
\newtheorem{xmp}[Theorem]{Example}
\newtheorem{cnj}[Theorem]{Conjecture}
\newtheorem{cst}[Theorem]{Construction}
\newtheorem{claim}{Claim}

\def \thm {\begin{Theorem} \ }
\def \ethm {\end{Theorem}}

\def \bp {\begin{prp} \ }
\def \ep {\end{prp}}

\def \bpm {\begin{prm} \ }
\def \epm {\end{prm}}

\def \bc {\begin{crl} \ }
\def \ec {\end{crl}}

\def \bl {\begin{lem} \ }
\def \el {\end{lem}}

\def \bd {\begin{defi} \ \rm }
\def \ed {\end{defi}}

\def \brm {\begin{rmk} \ }
\def \erm {\end{rmk}}

\def \bcj {\begin{cnj} \ }
\def \ecj {\end{cnj}}

\def \bcs {\begin{cst} \ }
\def \ecs {\end{cst}}

\def \bxm {\begin{xmp} \ \rm }
\def \exm {\end{xmp}}

\def \bcl {\begin{claim} \ }
\def \ecl {\end{claim}}

\title{Singular Ramsey and Tur\'an numbers}
\author{Yair Caro~\thanks{Deptartment of Mathematics,
 University of Haifa-Oranim, Tivon 36006, Israel.
  {\tt yacaro@kvgeva.org.il}}
 \qquad and \qquad
  Zsolt Tuza~\thanks{Alfr\'ed R\'enyi Institute of Mathematics,
        Hungarian Academy of Sciences, H--1053 Budapest,
 Re\'altanoda u.~13--15, Hungary; and
  Department of Computer Science and Systems Technology,
  University of Pannonia, 8200 Veszpr\'em, Egyetem u.~10,
 Hungary. {\tt tuza@dcs.uni-pannon.hu}.} }
\date{\small Latest update on \version}

\begin{document}

\maketitle

\begin{abstract}
We say that a subgraph $F$ of a graph $G$ is singular if the degrees $d_G(v)$ are
 all equal or all distinct for the vertices $v\in V(F)$.
The singular Ramsey number $\Rs(F)$ is the smalles positive integer $n$ such that,
 for every $m\ge n$, in every edge 2-coloring of $K_m$, at least one of the
 color classes contains $F$ as a singular subgraph.
In a similar flavor, the singular Tur\'an number $\Ts(n,F)$ is defined as the
 maximum number of edges in a graph of order $n$, which does not contain $F$
 as a singular subgraph.
In this paper we initiate the study of these extremal problems.
We develop methods to estimate $\Rs(F)$ and $\Ts(n,F)$, present tight
 asymptotic bounds and exact results.
\end{abstract}

\section{Introduction}\label{sec:intro}

In this paper we introduce a new type of Ramsey and Tur\'an numbers, where the classical condition
 of the occurrence of a specified subgraph in an edge-colored complete graph
 is combined with restrictions on vertex degrees in the monochromatic host graph.

\subsection{Brief survey on degree-constrained problems}

The smallest particular case of Ramsey's theorem is that on six vertices every graph or its
complement contains the triangle $K_3$. Starting from here,
Albertson \cite{Al94} proved\footnote{Considerable delay occurred between the birth and the
 publication of \cite{Al94}, and some of the follow-up papers appeared even several years earlier.}
 that for $n \ge 6$ in every 2-coloring of the edges of $K_n$ there is a monochromatic $K_3$ with
  two equal degrees.
Inspired by this result several papers were written, see for example
 \cite{Al97,AlBo94,BiWi05,BoCaFe,CaHaPe16}.

An obvious step after \cite{Al94} is to try to generalize this result to other graphs and also to try to bound the difference between the maximum and minimum degree of the specified monochromatic subgraph.
The efforts in the direction can be summarized as follows.

In \cite{AlBe91} Albertson and Berman showed that $K_n$ can always be colored red-blue in
 such a way that 
 no red $K_4$ occurs and no blue $K_2$ has equal monochromatic degree at its two vertices.
This shows that the phenomenon observed by Albertson is isolated and not extended to other graphs. But the authors of \cite{AlBe91} also showed that for $n \ge 6$ in every 2-coloring of $K_n$ there is
a $K_3$ with spread of the degrees at most 5, where the spread of a sequence
 $D = : \{ d_1,\dots,d_m \}$ is defined as 
  $\max \{  |d_i - d_ j| \mid 1 \le i , j \le  m\}$.
An extension of their result is presented in \cite{So97}.

In the papers \cite{ChSc93,ECRS93} Chen, Erd\H os, Rousseau and Schelp developed the notion of
 spread explicitly and proved in \cite{ECRS93}
  that every graph on at least $k+2$ vertices contains at least
 $k+2$ vertices whose degrees have spread at most $k$.
This is a non-trivial extension of the popular observation that every graph with
 more than one vertex has two vertices of the same degree.
From the quoted theorem the authors also proved among other things that
 for every graph $G$ and every $n \ge \R{G}$ (the classical Ramsey number)
 every 2-coloring of $K_n$ contains a monochromatic copy of $G$,
 whose vertex degrees in the host monochromatic graph have spread at most $\R{G} - 2$,
  and that in a certain sense this upper bound is tight.
An easy corollary is that the spread 5 from the Albertson--Berman result mentioned above
 can be reduced to 4 for $n \ge 6$, which is best possible (as already noted in \cite{ECRS93}).

Albertson  \cite{Al92} also introduced the corresponding Tur\'an number, namely
 the maximum number of edges in a graph on $n$ vertices having no copy of $K_m$
  with all degrees equal, and presented an exact bound.
(In an earlier paper \cite{CaErVi88} Caccetta, Erd\H os and Vijayan studied a Tur\'an-type problem
 concerning the existence of a complete graph $K_m$ with large degrees.)

A closely related subject is that of constant-degree independent sets, introduced by Albertson and Boutin \cite{AlBo94}, which was recently further developed by Caro, Hansberg and Pepper \cite{CaHaPe16}.
The latter considered various bounds on the constant-degree $k$-independent set in
trees, forest, $d$-degenerate graphs and $d$-trees.
Yet another direction concerns low-degree independent sets in planar graphs,
 developed by many authors and best presented in \cite{BiWi05}.

Further related notions are the so-called fair dominating sets (which actually are regular
 dominating sets, see Caro, Hansberg and Henning \cite{CaHaHe12}),
 irregular independence number and irregular domination number (Borg, Caro and Fenech \cite{BoCaFe}),
 and the problem of monochromatic degree-monotone paths in 2-colorings of the edges of complete
 graphs (Caro, Yuster and Zarb \cite{CaYuZa17}).

\subsection{Singular Ramsey and Tur\'an numbers}

Albertson and Berman \cite{AlBe91} presented edge 2-colorings of $K_n$ avoiding a monochromatic
 copy of $G$ with all monochromatic degrees in $K_n$ equal.
On the other hand, the opposite possibility of having a monochromatic copy of $G$
with all its vertices having distinct monochromatic degrees in $K_n$ is very easy to exclude,
 by any decomposition of $K_n$ into two regular spanning graphs $H$ and $\ov{H}$.
However, simultaneous exclusion of the two cases is impossible if $n$ is large.
This fact motivates our present study.

\bd
Let $k \ge 1$ be an integer. A sequence $a_1 \le a_2 \le \dots \le a_n$ of integers is
 called \emph{$k$-singular} if either $a_1 = \dots = a_n$ or  for every $j = 1,\dots,n-1$,
   $a_{j+1} - a_ j  \ge k$.
Also if $a_1,a_2,\dots,a_n$ are integers (repetitions are allowed),
 we say that they form a $k$-singular set
 if putting them in increasing order we obtain a $k$-singular sequence.
(Hence, ``set'' may mean ``multiset'' in this particular context.)
\ed

\bd
A subgraph $H$ of graph $G$ is called \emph{$k$-singular} if the degree sequence of its
 vertices in $G$ --- where $G$ is termed the \emph{host graph} --- forms a $k$-singular sequence.
\ed

For short, in case of $k=1$, a 1-singular sequence is called \emph{singular sequence},
 and a 1-singular subgraph is called \emph{singular subgraph}.

\bsk

Let now $\cF$ be a family of graphs.

\bd
The \emph{$k$-singular Ramsey number} $\Rs (\cF,k)$ is defined as the smallest integer $n$ such that
 in every 2-coloring of the edges of $K_m$ for any $m \ge n$,  one of the graphs induced by the
  color classes contains a $k$-singular member of $\cF$.
\ed

\brm
 If in a graph $G$ the subsequence of degrees belonging
to a set $B$ of vertices is $k$-singular, then so does the subsequence
 belonging to $B$ in the complement of $G$ as well. Hence,
 in case of two colors, the vertex sets of $k$-singular subgraphs in color 1
  coincide with those in color 2 (for any $k\ge 1$).
\erm

In a similar flavor, as a little deviation, we also introduce a Tur\'an-type function.

\bd
Given $\cF$, and a natural number $k$, the \emph{$k$-singular Tur\'an number}
 --- as a function of the order $n$ ---
 denoted by $\Ts(n,\cF,k)$ is defined as the maximum number of edges
  in a graph $G$ on $n$ vertices that contains no $k$-singular copy of any $F\in \cF$.
In particular, let $\Ts(n,q)$ be the maximum number of edges in a graph $G$  of order $n$
 that contains no singular copy of $K_q$.
\ed

For \emph{singular Ramsey numbers} we shall use the simpler notation $\Rs(\cF) = \Rs(\cF,1)$
 for $k=1$, and we write $\Rs(F)$ for $\Rs(\{F\})$.

It is also natural to introduce non-diagonal and multicolored versions of $\Rs(F)$.

\bd
If $F_1$ and $F_2$ are two graphs, their singular Ramsey number $\Rs(F_1,F_2)$ is
 the smallest $n$ such that, for every $m\ge n$, every 2-coloring of $K_m$ contains a
 singular copy of $F_1$ in the first color or a singular copy of $F_2$ in the second color.
More generally, also for an integer $s>2$, one may consider $s$ families $\cF_1,\dots,\cF_s$
 of graphs and define the $k$-singular Ramsey number $\Rs (\cF_1,\dots,\cF_s,k)$ as the
  smallest integer $n$ with the property that, for any $m \ge n$, in every coloring of the
   edges of $K_m$ with $s$ colors, there is an $i$ ($1\le i\le s$) such that the graph
    induced by the $i$th color class contains a $k$-singular\footnote{See the concluding
  section for some possible interpretations of this definition more precisely.}
     member of $\cF_i$.
\ed

\brm (Non-monotonicity.)
Let the number of colors be fixed.
If every coloring of $K_n$ contains a monochromatic singular copy of some $F\in\cF$,
 still there is no guarantee that so does every coloring of $K_{n+1}$ as well.
This issue concerning (non-) monotonicity was observed already in the first papers by
 Albertson, and ever since;
it is treated by imposing the condition for every  $m \ge n$ in the definition of $\Rs(\cF)$,
 rather than just taking the smallest $n$ forcing a singular monochromatic $F$ in every
 2-coloring of $K_n$.
\erm

\brm (Monotonicity Principle.)   \label{rm:mon}
 It is obvious --- but will be applied at some point below --- that the function $\Rs$ is
  monotone with respect to inclusion, for any fixed number of colors; for instance,
  if $F_1\ssq G_1$ and $F_2\ssq G_2$, then $\Rs(F_1,F_2)\leq \Rs(G_1,G_2)$ holds.
\erm

\brm
In the classical version of Ramsey and Tur\'an numbers, isolated vertices are
 practically irrelevant, namely $\R{G\cup mK_1}=\max(\R{G},m+|V(G)|)$;
 but this is not at all the case in the singular version.
For instance, it can easily be shown (partly following also from some later observations) that
 for the graph $G =  P_3 \cup  K_1$ --- the path on 3 vertices plus an isolated vertex ---
   we have $\Rs(P_3 \cup  K_1) = 10$ while $\R{P_3 \cup K_1} = 4$,
    moreover $\Rs(P_3)=5$ and $\R{P_3} = 3$.
   (Also, one may observe that $\Rs(3K_1)=5$ while $\R{3K_1}=3$.)
Similarly, the Tur\'an number of $K_2 \cup K_1$ is zero for every $n\ge 3$, but $K_4-e$
 does not contain it as a singular subgraph, therefore $\Ts(4,K_2 \cup K_1,1)=5$.
\erm

In this paper we will mostly consider Ramsey-type results for two colors,
 and develop a couple of methods
 suitable for determining the exact value of singular Ramsey numbers in both the diagonal and
 non-diagonal cases, provided that the specified graphs satisfy certain properties.
We also present asymptotic estimates, and the $k$-singular version will be touched, too.
In a section after the Ramsey-type results we provide tight asymptotics for the $k$-singular
 Tur\'an number of a graph.

\subsection{Our results}

While the star graphs can be considered as the easiest infinite class of graphs
 concerning the classical Ramsey numbers (they almost admit a one-line proof),
 they turn out to be a bit complicated in the singular version.
For this reason, although we present a complete solution, we do not discuss them
 earlier than in Section \ref{s:star}.
Before that, we give some general lower and upper bounds (Section \ref{s:bound}),
 describe some methods to derive tight estimates (Section~\ref{s:meth}),
 and determine exact results for all,
  but one, graphs with at most four vertices and edges,
 with the unique exception of $C_4$
 (Section \ref{s:small}).
Tight asymptotics for singular Tur\'an numbers are given in Section \ref{s:tur}.
Some open problems are mentioned in the concluding section.

\subsection{Terminology and notation}
   \label{ss:term-not}

\textit{Particular graphs.}
We use standard notation $P_n$ and $C_n$ for the path and the cycle on $n$ vertices;
 $K_{p,q}$ for the complete bipartite graph with $p$ and $q$ vertices in its classes;
 and $mK_2$ for the matching with $m$ edges.
The \textit{claw} is the graph $K_{1,3}$.
The \textit{paw}, which we abbreviate in formulas as $PW$, is the graph with four vertices
 and four edges obtained from $K_3$ by adding a pendant vertex (or from $K_4$ by
 deleting the edges of a $P_3$).
The \textit{bull} is the graph obtained from $K_3$ by adding two pendant vertices
 which are adjacent to two of its distinct vertices (a self-complementary graph with
 five vertices and five edges).

\textit{Vertex degrees.}
The degree of a vertex $v$ in a graph $G$ is denoted by $d_G(v)$, or simply $d(v)$
 if $G$ is clear from the context.
Minimum and maximum degree are denoted by $\delta(G)$ and $\Delta(G)$, respectively.
A \textit{degree class} consists of all vertices having the same degree; hence the
 degree classes partition $V(G)$, and their number is equal to the number of
 distinct values which occur in the degree sequence of $G$.
Given a vertex partition $V_1\cup\cdots\cup V_k=V(G)$, and a vertex $v\in V_i$,
 the \textit{internal degree} of $v$ is the number of its neighbors inside $V_i$,
 and its \textit{external degree} is the number of its neighbors in $V(G)\setminus V_i$.

\textit{Ramsey number.}
We denote the Ramsey number by $\R{\cF}$, that is the smallest $n$ such that
 in every 2-coloring of the edges of $K_n$, one of the color classes contains a
  monochromatic member of  $\cF$.

\textit{Substitution.}
Let $H$ be a graph with $k$ vertices $v_1,\dots,v_k$, and let $F_1,\dots,F_k$ be
 $k$ non-null graphs ($F_i=K_1$ is allowed).
The \textit{substitution} of $F_1,\dots,F_k$ into the ``host graph'' $H$,
 denoted by $H[F_1,\dots,F_k]$, is the
 graph whose vertex set is the disjoint union $V(F_1)\cup\dots\cup V(F_k)$, each
 $V(F_i)$ induces the graph $F_i$ itself, and two vertices $x\in V(F_i)$ and
 $y\in V(F_j)$ ($i\ne j$) are adjacent in $H[F_1,\dots,F_k]$ if and only if
 $v_iv_j$ is an edge in $H$.
In this construction we say that the graph $F_i$ is substituted for $v_i$.

 In a graph $G=(V,E)$, the subgraph induced by a set $Y\sst V$ is denoted by $G[Y]$.

\section{Singular Ramsey numbers: General bounds}
   \label{s:bound}

We start with the following easy lemma.

\bl   \label{l:1}
 Every sequence of $k(n-1)^2 +1$ integers contains a $k$-singular subsequence of cardinality at least $n$.
\el

\pf
 Suppose we have no $n$ equal elements in the sequence. Then we must have at least $k(n-1) + 1$ elements of distinct values.
Reorder them in increasing order, say $a_1 < \dots < a_ { k(n-1)+1}$.
Take the subsequence $a_{jk +1}$ for $j = 0,\dots,n-1$. Clearly  this is a $k$-singular $n$-term sequence.
\qed

\thm   \label{th:1}
For any two families $\cF_1,\cF_2$ of graphs and every natural number $k\ge 1$
 the following general upper bound holds:
  $$\Rs(\cF_1,\cF_2,k) \le k(\R{\cF_1,\cF_2} - 1)^2 +1.$$
\ethm

 \pf
Consider a 2-coloring of the edges of $K_m$, for any
 $m \geq k(\R{\cF_1,\cF_2} - 1) ^2 +1$.
Let $G_1$ and $G_2$ be the subgraphs obtained by the edges of color~1
 and color~2, respectively. 
By Lemma \ref{l:1} the sequence of degrees of the vertices of $G_1$
 contains a $k$-singular subsequence of cardinality $\R{\cF_1,\cF_2}$.
The degrees of the corresponding vertices form a $k$-singular subsequence
 also in $G_2$.
Now consider the 2-colring induced on the complete graph on those
 $\R{\cF_1,\cF_2}$ vertices.
By definition there is either a monochromatic copy of a graph $G \in \cF_1$
 in color 1 or of a graph $H \in \cF_2$ in color 2. 
Hence the degrees of $G$ (in the first case) form a $k$-singular subsequence in the
 host graph $G_1$ or the degrees of $H$ (in the second case) form a $k$-singular
  subsequence in the host graph $G_2$. 
Thus a required $k$-singular subgraph occurs whenever
 $m \geq k(\R{\cF_1,\cF_2} - 1) ^2 +1$, which means
 $\Rs(F_1,F_2,k) \leq k(\R{\cF_1,\cF_2} - 1) ^2 +1$.
 \qed

\bsk

An immediate corollary is:

\bc ~~ \label{c:square}
\begin{itemize}
\item[$(i)$] For every graph $G$ we have $\Rs(G) \le (\R{G} - 1)^2 +1$, and also
  $\Rs(G,H) \le (\R{G,H} - 1)^2 +1$ for any two graphs $G$ and $H$.

\item[$(ii)$] Every 2-coloring of $K_{k(n-1)^2+1}$ contains a monochromatic $k$-singular tree of order
 at least $n$.
\end{itemize}
\ec

\pf
$(i)$ This is just the case $\cF = \{G\}$, or $\cF_1 = \{G\}$ and $\cF_2 = \{H\}$,
  with $k = 1$ in Theorem \ref{th:1}.

$(ii)$ Consider the degree sequence in the graph induced by the edges colored  1.
 By Lemma \ref{l:1} there is a $k$-singular subsequence of $n$ degrees.
Consider now the induced coloring on the complete graph $K_n$  whose vertices are those forming the $k$-singular sequence.
Since every graph or its complement is connected, it follows that there is a connected monochromatic subgraph of order $n$ whose degree sequence is $k$-singular  in the host graph, and hence such a tree occurs.
\qed

\bsk

Having proved a general upper bound, we next supply a general quadratic lower bound.

\thm   \label{th:2}
Let $G$ be any graph on $n\ge 3$ vertices.
  Then $\Rs(G) \ge \max \{ \R{G} ,$ $ (n-1)^2 +1 \}$.
\ethm

\pf
Trivially $\Rs(G) \ge \R{G}$, so we only have to show $\Rs(G) \ge (n-1)^2 +1$.

We will construct a graph $H$ on $(n-1)^2$ vertices whose vertex set $V$ is partitioned into $n - 1$ subsets $V_0,\dots,V_{n-2}$, each of cardinality $n-1$, such that all vertices
 in $V_i$ have the same degree ($i = 0,1,\dots,n-2$) but vertices from distinct subsets have distinct degrees.
Then clearly no copy of $G$ in $H$ and $\ov{H}$ can be singular,  as it must take at least two vertices in the same class and at least two vertices in distinct classes.

If $n-1$ is even, then we simply insert any $i$-regular graph inside $V_i$.
 (Such graphs exist, e.g.\ by taking $i$ perfect matchings from
  any 1-factorization of $K_{n-1}$.)

If $n-1$ is odd, then depending on residue modulo 4, one of the sequences
 $0,1,\dots,n-2$ and $1,2,\dots,n-1$ contains an \textit{even} number of \textit{odd} terms.
If it is $0,1,\dots,n-2$, then we insert a regular graph of degree $2\cdot\lfloor \frac{i}{2} \rfloor$
 inside $V_i$ (e.g., the union of $\lfloor i/2 \rfloor$ edge-disjoint Hamiltonian
 cycles of $K_{n-1}$).
Moreover we insert a perfect matching between $V_1$ and $V_3$, between $V_5$ and $V_7$,
 ..., between $V_{n-5}$ and $V_{n-3}$.
Else, if it is $1,2,\dots,n-1$, then we insert a regular graph of degree $2\cdot\lfloor \frac{i+1}{2} \rfloor$
 inside $V_i$ (e.g., the union of $\lfloor \frac{i+1}{2} \rfloor$ edge-disjoint Hamiltonian
 cycles of $K_{n-1}$) and take
 a perfect matching between $V_0$ and $V_2$, between $V_4$ and $V_6$,
 ..., between $V_{n-4}$ and $V_{n-2}$.

These graphs satisfy the requirements, proving the lower bound for all $n$.
\qed

\brm
 An alternative proof --- which also works in the $k$-singular case for $k\ge 2$ ---
  can be obtained from the Erd\H os--Gallai characterization
  of graphical sequences.
 We note that for some combinations of $k$ and $n$ (both even)
  an analogous construction with $k(n-1)^2$ vertices is not possible,
  because a graph cannot have an odd number of odd-degree vertices.
\erm

In particular, the following bounds are obtained from the above estimates.

\bc   \label{naiv}
 If $G$ is a graph of order $n\ge 3$, then
  $$
    \max \{ \R{G} , (n-1)^2 +1 \} \le \Rs(G) \le (\R{G} - 1)^2 +1 .
  $$
\ec

If $\cG$ is a class of graphs in which $\R{G}$ is a linear function of $|V(G)|$ over all
 graphs $G\in\cG$, then the growth order of both estimates in Corollary \ref{naiv} is
 quadratic in $n$.
In particular, applying the theorem of \cite{CRSzT83} on the Ramsey numbers of graphs with bounded
 maximum degree, we obtain:

\thm   \label{th:4}
Let $\cG$ be the class of graphs with bounded degree $\Delta$ fixed.
Then for all $G\in\cG$ of order $n$ we have $\Rs(G)=\Theta(n^2)$, as $n\to\infty$.
\ethm

We conclude this section with a sufficient condition ensuring that the lower bound
 in Corollary \ref{naiv} holds with equality.
This result also exhibits a significant difference between the classical and the singular versions
 of Ramsey numbers concerning the role of isolated vertices.

\bp   \label{p:isol}
 Let $G=H\cup mK_1$, i.e.\ the graph obtained from a graph $H$ by adding $m$ isolated vertices.
 If $|V(G)|\ge \R{H}$,
  then $\Rs(G)=(|V(G)|-1)^2 + 1$.
\ep

\pf
 We only have to prove that $(|V(G)|-1)^2 + 1$ is an upper bound on $\Rs(G)$.
If $n\ge (|V(G)|-1)^2 + 1$, then in every 2-coloring of $K_n$ the subgraph of color 1 contains
 a singular subgraph, say $G^*$, on $|V(G)|=|V(H)|+m\ge \R{H}$ vertices.
Thus, a singular monochromatic copy of $H$ occurs, either in color~1 or in color 2,
 which can be supplemented to a singular
 copy of $G$ because the $m$ isolated vertices put no restriction on the color distribution
 in the rest of $G^*$.
\qed

\section{Some methods}
   \label{s:meth}

Assume that a graph $G$ has been fixed, for which we wish to find estimates on $\Rs(G)$.
We say that a graph $F$ is \textit{$G$-free} if $F$ does not contain any subgraph
 isomorphic to $G$.
Moreover, let us call $F$ an \textit{R-graph} for $G$ if both $F$ and
 $\overline{F}$ are $G$-free.
Analogously, we say that $F$ is an \textit{SR-graph} (`S' standing for `singular') for $G$
 if neither $F$ nor $\overline{F}$ contains a singular subgraph isomorphic to $G$.

Lower bounds on $\Rs(G)$ will be obtained by constructing SR-graphs from
 several (smaller) R-graphs.
We call this the technique of canonical colorings.
Possible different approaches will be described in the next two subsections,
  and a kind of combination of them afterwards.

The fourth subsection presents a method to derive upper bounds when some favorable
 information concerning the structure of R-graphs of order $\R{G}-1$ is available.
  This approach will lead to exact results in several cases.
Finally we mention another approach to upper bounds, based on vertex degrees.

\subsection{Non-regular Canonical Coloring, NRCC}

 This approach is useful when `large' R-graphs are not regular.
For instance, the claw $K_{1,3}$ and its complement $K_3\cup K_1$ are the two
 R-graphs of order 4 for $G=2K_2$, and also for $P_4$,
  but neither of them is regular.
We apply this method in Section \ref{ss:p4m2}.

Let $G$ be a graph on $n$ vertices, and let $t \le n-1$.
Consider $t$ copies of (not necessarily isomorphic)
 R-graphs over mutually disjoint vertex sets $V_1,\dots,V_t$.
Suppose that we can insert edges between the vertex classes (but not inside them)
 to obtain a graph $H$ with the following properties:
\begin{enumerate}
  \item In each vertex class $V_i$ ($i = 1,\dots,t$)  all the degrees $d_H(v)$ are equal.

  \item Degrees of vertices belonging to distinct vertex classes are distinct.
\end{enumerate}

\bl
With the assumptions above, we have $$\Rs(G) \ge |V(H)| + 1 = 1 + \sum_{i=1}^t  |V_i| . $$
\el

\pf
In such a case $H$ and $\ov{H}$ have exactly $t$ classes of distinct degrees and $t \le n-1$,
 hence no copy of $G$ with all degrees distinct is possible (there are too few distinct degree classes).
Also, since each set $V_i$ induces an R-graph in $H$,
 no copy of $G$ with all degrees equal is possible as it should be contained in a unique degree class.
Hence $H$ is an SR-graph, showing $\Rs(G) \ge |V(H)| + 1$.
\qed

\subsection{Regular Canonical Coloring, RCC}

We can apply this approach when there exist `large' R-graphs which are regular.
(The first classical example is $G=K_3$ whose unique largest R-graph is $C_5$.)
We apply this method in Section \ref{ss:triclaw}.

Let $F$ be an R-graph on $q$ vertices $v_1\dots,v_q$, and let $H_1,\dots,H_q$ be
 $q$ further R-graphs.
Denote by $H =F[H_1,\dots,H_q]$ the graph obtained by taking the vertex-disjoint
 copies of $H_1,\dots,H_q$ and making all the vertices of $H_i$ adjacent to
  all the vertices of $H_ j$ if and only if the vertices $v_i$ and $v_ j$ are adjacent in $F$.

Suppose that $H$ has the following properties:
\begin{enumerate}
  \item Each $H_i$ ($i = 1,\dots,q$) is a regular induced subgraph of $H$.

  \item If $i\neq j$, then the degrees $d_H(v)$ for vertices $v$
   in $H_i$ and $H_j$ are not the same.
\end{enumerate}

\bl
With the assumptions above, we have $$\Rs(G) \ge |V(H)| + 1 = 1 + \sum_{i=1}^q  |V(H_i)| .$$
\el

\pf
Observe first that since all vertices of $H_i$ are connected to the same vertices outside $H_i$
 and also have the same degree inside $H_i$, it follows that $H_i$ is a regular subgraph in $H$
  (hence the name Regular Canonical coloring).
Since $H_i$ and its complement $\ov{H_i}$ are $G$-free, property 2 implies
 that there is no copy of $G$ with all degrees equal.

If there was a copy of $G$ with all degrees distinct in $H$, then no $H_i$ would contain
 more than one vertex from $G$.
Hence, by the construction, there would be a copy of $G$ in $F$, but this is impossible because
 $F$ and $\ov{F}$ are $G$-free.

Thus $H$ is an SR-graph for $G$, and hence $\Rs(G) \ge |V(H)| +1$.
\qed

\subsection{A mixed construction}
   \label{ss:mixedconst}

 We apply this method in Sections \ref{ss:triclaw} and \ref{ss:paw}.

The construction starts with a graph $H$ such that $H$ contains no singular $G_1$
 and $\overline{H}$ contains no singular $G_2$.
Partition $V(H)$ into some number of subsets, say $V(H)=X_1\cup \cdots \cup X_k$;
 many of those $X_i$ may also be singletons.
As a generalization of substitution, we replace those $X_i$ with
 mutually vertex-disjoint graphs
 $Q_1,\dots,Q_k$ such that each $Q_i$ is regular, $G_1$-free, and $\overline{Q_i}$
 is $G_2$-free.
The plan is to create a graph $F$ whose degree classes are the sets $V(Q_i)$,
 using the structure of $H$.
If $X_i$ consists of $q_i$ vertices from $H$, then we partition $V(Q_i)$ into $q_i$
 subsets.
The vertices in the $j^\textrm{th}$ part of $Q_i$ are completely adjacent to those
 classes $Q_\ell$ which correspond to the neighbors of the $j^\textrm{th}$ vertex
 of $X_i$ in $H$.
(In particular, if $X_i$ and $X_\ell$ are singletons adjacent vertices, then
 we take complete bipartite adjacency between $Q_i$ and $Q_\ell$.)

A delicate detail in this approach is to ensure that two vertices have the same
 degree if and only if they are in the same $Q_i$.
This needs a careful choice of the orders $|V(Q_i)|$, the internal degree
 of each $Q_i$, and also the sizes of the partition classes inside $Q_i$.

\subsection{Ramsey-stable graphs}

 The tool described in this subsection will turn out to be substantial, in the
  proofs of upper bounds in several results below.

Let $G$ be a given graph for which we wish to determine or estimate the value of $\Rs(G)$.
Consider an R-graph $H$ for $G$, say with $k$ vertices $v_1,\dots,v_k$.
Let $N_i$ denote the set of vertices adjacent to $v_i$ (the neighborhood of $v_i$).

\bd
 We call $H$ a \textit{Ramsey-stable graph} for $G$ if, for each $1\le i\le k$,
  the unique way to obtain an R-graph of order $k$,
  in which $H-v_i$ is an induced subgraph, is to join a new vertex
  to all vertices of $N_i$, and not to join it to any other vertex of $H-v_i$.
 Ramsey-stable graphs for a pair $(G_1,G_2)$ of graphs can be defined analogously.
\ed

\bxm
 The 5-cycle is Ramsey-stable for $K_3$, and also for $K_{1,3}$, because the only way 
  to extend $P_4$ to an R-graph for $K_3$, or for $K_{1,3}$, is to
  join a new vertex to the two ends of $P_4$.
\exm

\brm   \label{r:regRstab}
 More generally than the previous example, if we know that all $n$-vertex R-graphs for
  a given $G$ are regular, then every R-graph $H$ of order $n$ is Ramsey-stable for $G$
   because exactly the vertices of minimum degree in $H-v$ have to be joined
   by an edge to the new vertex.
\erm

Assume that $F$ is an SR-graph for a given graph $G$, and that the degree sequence
 of $F$ contains precisely $k$ distinct values.
We partition $V(F)$ into the degree classes $V_1,\dots,V_k$.
Pick one (any) vertex $v_i$ from each class $V_i$, and denote by $H$ the graph induced by
 $\{v_1,\dots,v_k\}$ in $F$.
Since the set $\{v_1,\dots,v_k\}$ is irregular in $F$, we see that $H$ is an R-graph for $G$.

The significance of Ramsey-stable graphs is shown by the following lemma,
 which will be crucial in several proofs later on.
As a side product, it also implies that if a suitable choice of $\{v_1,\dots,v_k\}$ gives us
 a Ramsey-stable $H$, then all possible choices of the $v_i\in V_i$ ($i=1,\dots,k$)
  yield the same $H$.

\bl   \label{l:subst}
   (Regular Substitution Lemma.)
 Let $F,G,H$ be graphs as above.
   If $H$ is Ramsey-stable for $G$, then $F$ is obtained from $H$ by
     substituting a regular R-graph for each vertex $v_i$ of $H$.
 The same structure is valid when $H$ is Ramsey-stable for a pair $(G_1,G_2)$.
\el

\pf
 Assume that $H$ is Ramsey-stable for $G$; the case of $(G_1,G_2)$ can be handled
  in exactly the same way.
 Then for any $i$, replacing the vertex $V(H)\cap V_i$ with any $v\in V_i$, the neighborhood
  remains the same, by assumption.
 Hence every $v_j\in N_i$ (which has been taken from the degree class $V_j$)
  is completely adjacent to $V_i$.
 This is true also when we view the edge $v_iv_j$ from the other side, from $v_j$;
  therefore $v_i$ --- and each of its replacement vertices, $v\in V_i$ ---
  is adjacent to the entire $V_j$.
 Consequently, for each edge $v_iv_j$ of $H$, the edges between $V_i$ and $V_j$ in $F$
  form a complete bipartite graph spanning $V_i\cup V_j$.
 On the other hand, by the analogous argument for the non-edges of $H$, we see that
  if $v_iv_j$ is not an edge in $H$, then there are
  no edges between $V_i$ and $V_j$ in $F$.
 Thus, $F$ is generated by the operation of substitution.
 As a quantitative consequence, the external degrees of vertices in any one $V_i$ are all equal.

 Equal external degrees imply for a degree class that the internal degrees must also
  be equal.
 This implies regularity inside each $V_i$.
\qed

\subsection{Vertex degrees}

In some cases the following approach is useful in deriving upper bounds on $\Rs(G)$.
 We apply it in Section \ref{ss:triclaw}.

\bl
 If, for a given graph $G$, every SR-graph of order $n$ has minimum degree $\delta$,
  then there can be at most $n-2\delta$ degree classes.
\el

 \pf
Consider any SR-graph $F$ of order $n$, and
 let $k$ denote the number of its degree classes.
Then, concerning the minimum and maximum degree we have
 $$\delta+k-1 \le \delta(F)+k-1 \le \Delta(F) = n-1-\delta(\overline{F}) \le n-1-\delta$$
from where we obtain $k\le n-2\delta$.
 \qed

\bsk

Typically one can use this in the way that if $n$ is large then an SR-graph should
 have not only large minimum degree but also a large number of degree classes, from which
 a contradiction is derived to the above inequality, concluding that $\Rs(G)\le n$.

\section{Exact results on $\Rs (G)$ for small graphs}
   \label{s:small}

The smallest nontrivial cases are the path $P_3=K_{1,2}$ and its subgraphs;
 they allow a simple solution for $k$-singular Ramsey numbers for all $k$,
  which we present in the first subsection.
In this way $K_3$ remains the unique graph $G$ of order three for which we
 do not know $\Rs(G,k)$ over the entire range of $k$.

All other subsections of this section deal with the case $k=1$ for small graphs, determining
 $\Rs$ for every graph with at most four vertices and at most four edges,
  except for $C_4$ where we have a non-trivial lower bound.
This also includes small star graphs (the claw $K_{1,3}$, and the $K_{1,2}$ which is
 treated under the name $P_3$); a general theorem for stars will be presented
 in Section \ref{s:star}.

\subsection{The path $P_3$ for general $k$ of singularity}

\thm   \label{th:5}
$\Rs(3K_1 , k) = \Rs(K_2 \cup K_1 , k) = \Rs(P_3 , k) = 4k +1$.
\ethm

\pf
Clearly, by Theorem \ref{th:1} above we get $\Rs(P_3 , k) \le 4k +1$ as $\R{P_3} = 3$.

For the lower bound consider the graph $H(k)$ on $4k$ vertices defined as follows:
$V(H(k))  = A \cup B$, where $A = \{ a_1,\dots,a_{2k} \}$, $B = \{ b_1,\dots,b_{2k}\}$,
 and $a_i$ is adjacent to $b_ j$  precisely when $i \le j $.

In this graph,
  which treats the lower bound for the three graphs $3K_1$, $K_2 \cup K_1$, $P_3$ together,
 every degree between 1 and $2k$ is repeated exactly twice,
 i.e.\ no triple can have equal degrees.
Also there cannot occur any $k$-singular subgraph of order three, because this would require
 that $\Delta(H(k)) - \delta(H(k)) \ge 2k $, however in $H(k)$ and hence also
  in its complement the difference is just $2k-1$.
\qed

\subsection{The path $P_4$ and the 2-matching $2K_2$}
   \label{ss:p4m2}

Here we prove:

\thm   \label{th:6}
$\Rs(2K_2) = \Rs(2K_2,P_4) = \Rs(P_4) = 13$.
\ethm

 \pf
By the Monotonicity Principle we have $\Rs(2K_2) \le \Rs(2K_2,P_4) \le \Rs(P_4)$,
 therefore it suffices to prove that $\Rs(2K_2)\ge 13$ and $\Rs(P_4)\le 13$.

For the lower bound on $\Rs(2K_2)$ we construct an SR-graph on 12 vertices.
Consider $V_1,V_2,V_3$, where $|V_i| = 4$  for $i=1,2,3$.
Let each of
$V_1 ,V_2  ,V_3$ induce a $K_3$ with an isolated vertex. The vertices are labeled as $V_1 = \{ x_1,x_2,x_3, x\}$ where $x$ is the  isolated vertex not in the $K_3 $, similarly $V_2  = \{ y_1,y_2,y_3 ,y\}$ with $y$ not in the $K_3 $, and $V_3 = \{ z_1,z_2,z_3 ,z \}$ with $z$ not in the $K_3$.

We complete these vertex classes to a graph $G$ (color 1) such that all degrees in $V_1$ are 7, all degrees in $V_2$ are 5, and all degrees in $V_3$ are 4.
Once this shall be done, there will be no copy of $2K_2$ with all degrees  equal in $G$ and
 neither in $\ov{G}$ because each $V_i$ induces $K_3 \cup  K_1$ in $G$ and $K_{1,3}$ in $\ov{G}$.
Also there will be no $2K_2$ with all degrees distinct since this would require four different degrees,
 while in both $G$ and $\ov{G}$ there are only three.
 We shall do the construction step by step.

First, connect $x_1$ to $y_1,y_2, y$; $x_2$ to $y_2,y_3,y$; $x_3$ to $y_1,y_3,y$;  and $x$ to $y_1,y_2,y_3,y$.
 The degrees are now 4 for $x,y$; 5 for $x_1,x_2,x_3,y_1,y_2,y_3$; 0 for $z$; and 2 for $z_1,z_2,z_3$.
Next,
connect $x_1$ to $z_1 ,z$; $x_2$ to $z_2 ,z$; $x_3$ to $z_3 ,z$; $x$ to $z_1,z_2 ,z_3$.
 Then the degrees are 7 for $x, x_1,x_2,x_3$; 5 for $y_1,y_2,y_3$; 4 for $y,z_1,z_2,z_3$; and 3 for $z$.
Finally,
connect $y$ to $z$, and we are done.

 For the upper bound on $\Rs(P_4)$ we will apply the Regular Substitution Lemma.
On four vertices precisely two graphs are R-graphs for $P_4$: the claw $K_{1,3}$ and its complement, the
 triangle $K_3$ with an isolated vertex.
  On five vertices every edge 2-coloring contains a monochromatic $P_4$.
Observe that each of $K_{1,3}$ and $K_3\cup K_1$ is a Ramsey-stable graph for $\Rs(P_4)$,
 because a 3-vertex subgraph with zero or two edges is extendable only to the claw, whereas
 that with one or three edges is extendable only to the triangle;
  either extension is unique
 also concerning the set of neighbors of the new vertex.

Suppose now for a contradiction that there exists an SR-graph $F$ for $P_4$ on
 at least 13 vertices.
There can be at most four degree classes in $F$, each on at most four vertices.
It follows that there are precisely four vertex classes.
Due to Lemma \ref{l:subst}, each degree class should induce a regular R-graph;
 but this is impossible for a class with four vertices, which must occur if $|V(F)|>12$.
This contradiction completes the proof.
 \qed

\subsection{The triangle $K_3$ and the claw $K_{1,3}$}
   \label{ss:triclaw}

Although there is no containment relation between $K_3$ and $K_{1,3}$, the
 unique R-graph of order 5 for both of them is the 5-cycle.
Moreover, on four vertices, every R-graph has positive minimum degree.
These facts allow us to treat the two graphs together,
 and prove the following theorem.

\thm   \label{th:7}
 $\Rs(K_3) = \Rs(K_{1,3}) = 22$.
\ethm

\nin
{\bf Proof of Lower Bound 22.}
 We construct an SR-graph of order 21.
Consider the 5-cycle $v_1v_2v_3v_4v_5$ as host graph, and substitute $F_1,\dots,F_5$
 for $v_1,\dots,v_5$ as follows:
  $$
    F_1 \cong F_2 \cong F_3 \cong C_5 \ , \quad F_4 \cong 2K_2 \ , \quad F_5 \cong K_2 \ .
  $$
Then the degrees are:
 \begin{itemize}
   \item for $F_1$, internal: 2, external: $2+5=7$, total: 9;
   \item for $F_2$, internal: 2, external: $5+5=10$, total: 12;
   \item for $F_3$, internal: 2, external: $5+4=9$, total: 11;
   \item for $F_4$, internal: 1, external: $5+2=7$, total: 8;
   \item for $F_5$, internal: 1, external: $4+5=9$, total: 10.
 \end{itemize}
Since the host graph and also the subgraphs substituted for the degree classes are
 $K_3$-free and $K_{1,3}$-free, no singular $K_3$ or $K_{1,3}$ occurs.
  \qed
\bsk

\nin
 {\bf Proof of Upper Bound 22.}
Since $\R{K_3} = \R{K_{1,3}} = 6$, we infer from Theorem \ref{th:1} that
 $\Rs(K_3) \le 26$ as well as $\Rs(K_{1,3}) \le  26$.
So we have to cover the cases $n = 22 ,23 ,24 ,25$, to show that a singular triangle
 and a singular claw necessarily occurs in each case.

For a contradiction, consider an SR-graph; we know that it can have at most five
 degree classes, each with at most five vertices.
Hence, the following combinations might occur:
 \begin{itemize}
   \item $22 = 5+5+5+5+2$
   \item $22 = 5+5+5+4+3$
   \item $22 = 5+5+4+4+4$
   \item $23 = 5+5+5+5+3$
   \item $23 = 5+5+5+4+4$
   \item $24 = 5+5+5+5+4$
   \item $25 = 5+5+5+5+5$
 \end{itemize}
We will show that all of them are impossible.

\bsk

\nin
 {\bf First Proof --- Degree Counting.}
We arrange the degree classes in decreasing order of size $|V_1|\ge|V_2|\ge|V_3|\ge|V_4|\ge|V_5|$,
 and denote the vertices of $V_i$ as $v_{i1}$ $v_{i2}$, $\dots .$\
Then consider the seven cases separately.
 \begin{itemize}
   \item $22 = 5+5+5+5+2$ --- Vertex $v_{11}$ has precisely two neighbors in each of the sets
     $\{v_{12},v_{13},v_{14},v_{15}\}$, $\{v_{21},v_{31},v_{41},v_{51}\}$,
      $\{v_{22},v_{32},v_{42},v_{52}\}$; and has at least one neighbor in each of
       $\{v_{23},v_{33},v_{43}\}$, $\{v_{24},v_{34},v_{44}\}$, $\{v_{25},v_{35},v_{45}\}$.
     Thus $d(v_{11})\ge 9$.
     Similarly, $v_{51}$ has precisely two neighbors in each of the sets
      $\{v_{1j},v_{2j},v_{3j},v_{4j}\}$ for $j=1,\dots,5$, hence $d(v_{51})\ge 10$.
     This means $\delta(F)\ge 9$, as the positions of the other vertices are analogous;
      and since we have five degree classes, $\Delta(F)\ge 13$ follows.
     The same inequalities must hold for $\overline{F}$, too.
     But $\Delta(F)\ge 13$ implies $\delta(\overline{F}) = |V(F)|-1-\Delta(F)\le 8$,
      a contradiction.
   \item $22 = 5+5+5+4+3$ --- Here $v_{11}$ has two neighbors in $\{v_{12},v_{13},v_{14},v_{15}\}$
      and also in each of the sets $\{v_{2j},v_{3j},v_{4j},v_{5j}\}$ for $j=1,2,3$;
      and has at least one neighbor in $\{v_{24},v_{34},v_{44}\}$, which means $d(v_{11})\ge 9$.
     Vertex $v_{41}$ has two neighbors in $\{v_{1j},v_{2j},v_{3j},v_{5j}\}$ for $j=1,2,3$,
      and at least one neighbor in each of $\{v_{14},v_{24},v_{34}\}$, $\{v_{15},v_{25},v_{35}\}$,
       $\{v_{42},v_{43},v_{44}\}$.
     Vertex $v_{51}$ has two neighbors in $\{v_{1j},v_{2j},v_{3j},v_{4j}\}$ for $j=1,2,3,4$,
      and at least one neighbor in $\{v_{15},v_{25},v_{35}\}$.
     Thus, $\delta(F)\ge 9$, a contradiction again.
   \item $22 = 5+5+4+4+4$ --- Here the vertices of $V_1\cup V_2$ must have degree
      at least 10, and the vertices of $V_3\cup V_4\cup V_5$ must have degree at least~9.
   \item $23 = 5+5+5+5+3$ --- Here the vertices of $V_5$ have two neighbors in
      $\{v_{1j},v_{2j},v_{3j},v_{4j}\}$ for $j=1,2,3,4,5$, while the other vertices have
       two neighbors in each of four 4-tuples and one neighbor in each of two triples.
     Thus $\delta(F)\ge 10$, $\Delta(F)\ge 14$, $\delta(\overline{F}) < 10$ ---
      a contradiction.
   \item $23 = 5+5+5+4+4$ --- Also here, every vertex has degree at least 10, hence the
      maximum degree should be at least 14.
   \item $24 = 5+5+5+5+4$ --- Here $\delta(F)\ge 11$ and $\Delta(F)\ge 15$ should hold.
   \item $25 = 5+5+5+5+5$ --- This graph should be 12-regular, despite that it has
      five degree classes. \qed
 \end{itemize}

\bsk

\nin
 {\bf Second Proof --- Ramsey-Stable Graphs.}
Since the 5-cycle is the unique R-graph on five vertices, any 5-tuple with one vertex
 from each degree class must induce $C_5$.
Thus, by the Regular Substitution Lemma, $F$ is obtained by substituting \textit{regular}
 R-graphs into $C_5$.
In particular, each 5-element $V_i$ must induce the 2-regular $C_5$.

The partition $22 = 5+5+5+5+2$ cannot occur because vertices in both neighbors of the 2-element class
 along the 5-cycle have external degree $5+2=7$ and internal degree 2, contradicting that they are
 distinct degree classes.
The same argument excludes $23 = 5+5+5+5+3$, $24 = 5+5+5+5+4$, and of course $25 = 5+5+5+5+5$ as well.

For $22 = 5+5+4+4+4$ note further that a 4-element $V_i$ must induce the 2-regular $C_4$ or the
 1-regular $2K_2$.
Thus, all internal degrees are between 1 and 2, and all external degrees are between 8 and 10,
 leaving room for no more than four degree classes while  we should have five of them.
For this reason, the case $22 = 5+5+4+4+4$ cannot occur, and the same argument
 excludes $23 = 5+5+5+4+4$.

The only case that remains is $22 = 5+5+5+4+3$.
External degree 7 can only occur for a 5-element class which has internal degree 2.
External degree 8 can only occur for a 4-element or a 5-element
 class, both having internal degree at least 1.
All other possibilities yield external degrees at least 9, thus $\delta(F)\ge 9$,
 and we can conclude as in the first proof that $\Delta(F)\ge 13$ should hold, from which we arrive
 at the contradiction $\delta(\overline{F})\le 8$.
\qed

\bsk

It turns out that the non-diagonal singular Ramsey number
 $\Rs(K_3,K_{1,3})$ is bigger.

\thm   \label{t:tricla}
 $\Rs(K_3,K_{1,3}) = 29$.
\ethm

\nin
 {\bf Proof of Lower Bound 29.}

We construct a graph $F$ on 28 vertices, without singular triangles, whose complement
 $\overline{F}$ does not contain any singular claws.
This $F$ will have $k=5$ degree classes $V_1,\dots,V_5$, where $|V_1|=|V_3|=|V_5|=6$
 and each of those classes induces $K_{3,3}$, while $|V_2|=|V_4|=5$ and both classes
 induce $C_5$.
Hence each degree class is internally regular, with no $K_3$ in it, and no $K_{1,3}$
 in the complementary graph.

We also partition $V_5$ into two sets as $V_5=V'\cup V''$, with the only restriction
 that $|V'|=4$ and $|V''|=2$, but no condition on the actual position of vertices.
The other edges of $F$ establish complete adjacencies
 \begin{itemize}
  \item between $V_1\cup V_2$ and  $V_3\cup V_4$,
  \item between $V_1\cup V_2$ and  $V'$,
  \item between $V''$ and  $V_3\cup V_4$.
 \end{itemize}
There are no other edges in $F$. Then the degrees are:
 \begin{itemize}
  \item in $V_1$: internal 3, external $6+5+4=15$, total 18;
  \item in $V_2$: internal 2, external $6+5+4=15$, total 17;
  \item in $V_3$: internal 3, external $6+5+2=13$, total 16;
  \item in $V_4$: internal 2, external $6+5+2=13$, total 15;
  \item in $V_5$, for vertices in any of $V'$ and $V''$: internal 3, external $6+5=11$, total 14.
 \end{itemize}
One can observe that every triangle of $F$ contains two vertices in the same $V_i$
 and one vertex in another class, hence no singular triangles occur.
Similarly the complement of $F$ contains no singular claws.
 Thus $F$ satisfies all requirements and yields $\Rs(K_3,K_{1,3}) \ge 29$.
\qed

\bsk

Concerning the upper bound we first observe some structural properties of the graphs
 which are R-graphs for $(K_3,K_{1,3})$.

\bcl   \label{cl:6}
 We have $\ndR(K_3,K_{1,3})=7$, and the unique R-graph of order 6 is $K_{3,3}$.
\ecl

 \pf
Observe that $K_{3,3}$ is the unique triangle-free graph of order 6 whose
 minimum degree is at least 3.
On the other hand, if the minimum degree is smaller than 3, then
 the complement contains $K_{1,3}$.
 \qed

\bcl   \label{cl:5}
 On five vertices there are precisely two graphs $H$ --- namely $C_5$ and $K_{2,3}$ ---
  such that $H$ is triangle-free and $\overline{H}$ is $K_{1,3}$-free.
 The first one, $C_5$, is a Ramsey-stable graph for $(K_3,K_{1,3})$.
\ecl

 \pf
All vertex degrees must be at least 2 (otherwise $\overline{H}$ contains $K_{1,3}$)
 and at most 3 (otherwise $H$ contains $K_3$ or $\overline{H}$ contains $K_{1,3}$).
If $H$ is 2-regular, then $H\cong C_5$.
In the remaining case assume that $d(v)=3$.
The three neighbors of $v$ must be mutually non-adjacent, otherwise $K_3\sst H$;
 and all of them have to be adjacent to the fifth vertex, since $\delta(H)\ge 2$.
  No further edges can occur, hence $H\cong K_{2,3}$ in this case.
Since no other R-graphs are possible, and $P_4$ is not an induced subgraph of
 $K_{2,3}$, it is clear that $C_5$ is Ramsey-stable.
 \qed

\bcl   \label{cl:4}
 Among the \underline{regular} four-vertex graphs $H$ there are precisely two --- namely $C_4$
  and $2K_2$ --- such that $H$ is triangle-free and $\overline{H}$ is $K_{1,3}$-free.
\ecl

 \pf
The other two regular graphs of order 4 are $K_4$ which contains $K_3$, and $4K_1$ whose
 complement contains $K_{1,3}$.
 \qed

\bsk

\nin
 {\bf Proof of Upper Bound 29.}

Let $F$ be an SR-graph for $(K_3,K_{1,3})$, say on $n:=\Rs(K_3,K_{1,3})-1$ vertices;
 we need to prove that $n\le 28$.
We see from Claim \ref{cl:6} that $F$ has at most six degree classes $V_1, \dots, V_k$,
 and $|V_i|\le 6$ holds for each of them.

\msk

\nin
 {\bf Case 1: Four vertex classes.}

This case is obvious: since $|V_i|\le 6$ holds for all $i$, we cannot have
 more than 24 vertices.

\msk

\nin
 {\bf Case 2: Six vertex classes.}

Picking one vertex $v_i$ from each vertex class $V_i$ we obtain an R-graph $H$
 of order 6.
Due to Claim \ref{cl:6}, we have $H\cong K_{3,3}$, which is Ramsey-stable.
The Regular Substitution Lemma implies that $F$ is obtained
 by substituting regular R-graphs for the vertices of $H$;
 the possible subgraphs with more than three vertices are listed in Claims
 \ref{cl:6}, \ref{cl:5}, and \ref{cl:4} (and $K_{2,3}$ is excluded).
Let us denote the subgraphs substituted into the partite sets by $Q_1,Q_2,Q_3$
 and $R_1,R_2,R_3$; and let their respective orders be $q_1,q_2,q_3,r_1,r_2,r_3$.
Also let us write $d(q_1),d(q_2),d(q_3),d(r_1),d(r_2),d(r_3)$ for their internal degrees.
We fix an indexing such that $d(q_1)\ge d(q_2)\ge d(q_3)$ and $d(r_1)\ge d(r_2)\ge d(r_3)$.
Note that all these $d$ are between 0 and 3 (and 0 can occur only if the substituted
 graph has at most three vertices).

Denoting $q=q_1+q_2+q_3$ and $r=r_1+r_2+r_3$, the degree set of $F$ is
 $$
   q+d(r_1) , \quad q+d(r_2) , \quad q+d(r_3) , \quad r+d(q_1) , \quad r+d(q_2) , \quad r+d(q_3) 
 $$
  with six mututally distinct values.
In particular, we must have strict inequalities $d(q_1)>d(q_2)>d(q_3)$ and
 $d(r_1)>d(r_2)>d(r_3)$.
It follows on each side of $K_{3,3}$ that each of $K_{3,3}$ and $C_5$ can be
 substituted only once, which implies $\max(q,r)\le 15$.
Moreover, assuming $q\ge r$ the degrees cannot be smaller than $r$ and cannot be
 larger than $q+3$, hence the presence of six distinct degrees implies $q+3\ge r+5$,
  i.e.\ $r\le q-2\le 13$.
Thus $n=p+r\le 28$.

\msk

\nin
 {\bf Case 3: Five vertex classes.}

As above, we pick one (any) vertex $v_i$ from each $V_i$ ($1\le i\le 5$),
 and consider the graph $H$ induced by them in $F$.
  Due to Claim \ref{cl:5}, this $H$ must be $C_5$ or $K_{2,3}$.
  Since $C_5$ is Ramsey-stable, the proof for it is easy.
Indeed, as above, the Regular Substitution Lemma implies that $F$ is
 obtained by substituting regular R-graphs for the vertices.
But $n=29$ or $n=30$ would imply that along the 5-cycle four
 consecutive substitutions would be $K_{3,3}$.
The two middle ones of them would have external degree 12, internal degree 3,
 total degree 15, contradicting the assumption that they form distinct
 degree classes.

 Hence, from now on we assume that $H\cong K_{2,3}$.
Re-label the indices, if necessary, so that the 2-element class of $K_{2,3}$ is
 $\{v_1,v_2\}$ and the 3-element class is $\{v_3,v_4,v_5\}$.
Although $K_{2,3}$ is not Ramsey-stable, vertices $v_1$ and $v_2$ have the
 property that replacing any one of them with a vertex from its class,
 we must obtain again a $K_{2,3}$, which implies that $\{v_3,v_4,v_5\}$ is
  completely adjacent to $V_1\cup V_2$.
Due to the exclusion of singular $K_3$,
 this also forces that $V_1$ and $V_2$ are completely non-adjacent.

For a vertex $v\in V_i$ from $V_3\cup V_4\cup V_5$ there can occur two situations:
 $v$ is adjacent either to $v_1$ and $v_2$ --- in which case we say that
 $v$ is in a \textit{stable position} --- or to the two vertices of
 $\{v_3,v_4,v_5\}\smin\{v_i\}$.

If all $v \in V_3\cup V_4\cup V_5$ are in a stable position, then we have
 complete adjacency between $V_1\cup V_2$ and $V_3\cup V_4\cup V_5$, moreover
 no edges can occur between $V_3$ and $V_4$, between $V_3$ and $V_5$, and
  between $V_4$ and $V_5$, and also between $V_1$ and $V_2$ either.
 This yields regular external degrees for each $V_i$.
Hence the internal degrees of $V_3$, $V_4$, and $V_5$ must be regular and
 mutually distinct, as well as those in $V_1$ and $V_2$,
  what implies that $|V_1\cup V_2|\le 6+5=11$ and
  $|V_3\cup V_4\cup V_5|\le 6+5+4=15$, thus $n\le 26$.

The occurrance of vertices in non-stable position requires a little more
 structural analysis.
For this, suppose that a $v\in V_5$ is adjacent to $v_3$ and $v_4$,
 instead of $v_1$ and $v_2$.
Since $\overline{F}$ contains no singular claws, and $v\in V_5$ already has
 non-neighbors in $V_1$ and $V_2$, all vertices of $V_3\cup V_4$ are
 adjacent to $v$.
Then no edges can occur between $V_3$ and $V_4$, for otherwise $F$ would
 contain a singular $K_3$.
Similarly, $v$ has no neighbors in $V_1\cup V_2$, because such a neighbor
 and $v$ would form a singular $K_3$ with $v_3$ (and also with $v_4$).

\vbox{
The non-adjacency of $V_3$ and $V_4$ also implies that all vertices in
 $V_3\cup V_4$ are in a stable position.
Thus, we have the following structure:
}

 \begin{itemize}
   \item there is complete adjacency between $V_1\cup V_2$ and $V_3\cup V_4$;
   \item $V_5$ admits a partition $V'\cup V''$ such that
     $V_1\cup V_2$ is completely adjacent to $V'$ and
     $V_3\cup V_4$ is completely adjacent to $V''$;
   \item no other edges occur between any $V_i$ and $V_j$ for $1\le i<j\le 5$.
 \end{itemize}
This structure implies that vertices in $V_1$ and $V_2$ have the same external degree,
 namely $|V_3|+|V_4|+|V'|$; and similarly, both $V_3$ and $V_4$ have external degree
  $|V_1|+|V_2|+|V''|$.
As a consequence, $|V_1|+|V_2|\le 6+5=11$ and $|V_3|+|V_4|\le 6+5=11$,
 and finally $n\le 22+|V_5|\le 28$.
\qed

\subsection{The paw graph}
   \label{ss:paw}

As another small graph, we determine the singular Ramsey number of the paw,
 that is a triangle with a pendant edge.
Its Ramsey number is $\R{PW}=7$.
Let us first summarize some facts about the R-graphs.

\bl   \label{l:claw}
 For the paw graph,
 \begin{itemize}
   \item[$(i)$] every graph on at most three vertices is an R-graph, and among them,
   the regular ones are $K_3$ and ist complement;
   \item[$(ii)$] on four vertices there are two regular R-graphs, the 1-regular $2K_2$
   and the 2-regular $C_4$;
   \item[$(iii)$] on five vertices there are three R-graphs, namely $K_2\cup K_3$ and its
   complement $K_{2,3}$ which are non-regular, and the 2-regular $C_5$;
   \item[$(iv)$] on six vertices there are two R-graphs, the 2-regular $2K_3$
   and its 3-regular complement, $K_{3,3}$.
 \end{itemize}
\el

 \pf
Parts $(i)$ and $(ii)$ are obvious.
For $(iii)$ one may note that $C_5$ is the unique R-graph for $K_3$ on five vertices, and
 of course it is an R-graph for the paw, too.
If $G$ contains a triangle and is an R-graph for the paw, then the triangle is a
 connected component.
This implies that on five vertices the complement of $G$ must contain $K_{2,3}$, and
 on six vertices the complement must contain $K_{3,3}$.
Then there cannot be any further edges in $\overline{G}$, hence $G \cong K_2\cup K_3$
 or $G \cong 2K_3$.
Analogously, if a triangle occurs in $\overline{G}$, then
 $G \cong K_{2,3}$ or $G \cong K_{3,3}$.
 \qed

\bsk

The quadratic formula yields the upper bound 37 on $\Rs(PW)$,
 but in fact the exact value is much smaller.

\thm
 $\Rs(PW)=31$.
\ethm

\nin
 {\bf Proof of Lower Bound 31.}

We construct a graph $F$ of order 30 which is an SR-graph for the paw.
 It will have five degree classes $V_1,\dots,V_5$, each of cardinality 6.
The degree classes induce R-graphs: $F[V_1] \cong F[V_3] \cong F[V_5] \cong 2K_3$,
 and $F[V_2] \cong F[V_4] \cong K_{3,3}$.
  (One may verify in the proof below that it would be equally fine
   to take $F[V_5] \cong K_{3,3}$.)
Further, we partition $V_5$ as $V_5 = V' \cup V''$, with $|V'|=2$ and $|V''|=4$.

We make complete adjacencies between any two of the three sets $V_1, V_2, V'$;
 and also between any two of $V_3, V_4, V''$.
There are no further adjacencies; i.e., the only edges between
 $V_1\cup V_2\cup V'$ and $V_3\cup V_4\cup V''$ occur inside $V_5$
 (namely between $V'$ and $V''$).

This $F$ contains no regular paw, because the degree classes are paw-free;
 and it has no irregular paw either, because omitting the internal edges of the
 degree classes (which edges certainly cannot occur in \textit{any} irregular subgraph)
 we obtain a graph which is generated by substituting independent sets into $2K_3$.
Now we have the following degrees:
 \begin{itemize}
   \item in $V_1$: external $6+2=8$, internal 2, total 10;
   \item in $V_2$: external $6+2=8$, internal 3, total 11;
   \item in $V_3$: external $6+4=10$, internal 2, total 12;
   \item in $V_4$: external $6+4=10$, internal 3, total 13;
   \item in $V_5$: external $6+6=12$, internal 2, total 14.
 \end{itemize}
Hence, $F$ satisfies all requirements and yields $\Rs(PW)>30$.
 \qed

\bsk

\nin
 {\bf Proof of Upper Bound 31.}
Suppose for a contradiction that there exists an
 SR-graph $F$ on at least 31 vetices.
We know that each degree class contains at most six vertices, therefore we have
 exactly six degree classes $V_1,\dots,V_6$.
Picking one vertex $v_i$ from each $V_i$, we get an R-graph, say $H$, of order 6, which must be
 either $2K_3$ or $K_{3,3}$, due to Lemma \ref{l:claw}$(iv)$.
Turning to the complement of $F$ if necessary, we may assume without loss of generality
 that $H\cong 2K_3$.

Since $2K_3$ is Ramsey-stable, we see from the Regular Substitution Lemma that
 $F$ is obtained by substiting regular R-graphs for the vertices of $H$.
We are going to analyze the feasible substitutions
 which create three distinct degree classes for each of the two components.
We shall see that it is not possible to have more than 16 vertices in a component,
 there is just one way to obtain 16, and there are exactly two ways to obtain 15.
Hence only the combinations $32=16+16$ and $31=16+15$ would yield $n>30$, but the
 argument below will show that each of them would force equal degrees to at least
 two of the $V_i$, which contradicts the definition of degree class.

\msk

\underline{18:} The unique feasible partition is $18 = 6+6+6$.
But then two of the degree classes induce the same R-graph ($2K_3$ or $K_{3,3}$),
 therefore they have the same degree in $F$, a contradiction.

\msk

\underline{17:} The unique feasible partition is $17 = 6+6+5$.
Then the vertices in both 6-classes have external degree 11.
In order that they have different degrees in $F$, one of them must induce $K_{3,3}$
 and the other induce $2K_3$.
Then their degrees in $F$ are 14 and 13, respectively.
However, the class of five vertices has external degree 12 and internal degree 2, yielding total 14,
 which is not feasible.

\msk

\underline{16:} The two feasible partitions are $16=6+6+4$ and $16=6+5+5$.
The latter is easy to exclude, because both 5-classes have internal degree 2 and external degree 11.
Concerning $16=6+6+4$ we see that the two subgraphs for `6' must have distinct internal degrees,
 hence one of them is $2K_3$, the other is $K_{3,3}$.
Both have external degree $6+4=10$, hence the vertex degrees are 12 and 13, respectively.
This implies that the subgraph for `4', which has external degree 12, must be $C_4$
 because $2K_2$ with internal degree 1 would repeat the degree 13.
Thus the degrees necessarily are
 $$
   12, \ 13, \ 14.
 $$
Of course, this cannot occur on more than one triangle; i.e., the case $n=32=16+16$ is impossible.

\msk

\underline{15:} The possible partitions are $15=6+6+3$, $15=6+5+4$, $15=5+5+5$.
We can immediately exclude the last one because `5' necessarily means $C_5$ with internal
 degree 2, hence in a substitution of the type $5+5+5$ the graph would be regular of degree 12.
Concerning $15=6+6+3$ --- similarly to the case of $16=6+6+4$ --- we see that $2K_3$ and $K_{3,3}$
 have to be substituted for $6+6$, yielding vertex degrees 11 and 12.
For `3' we have external degree 12, hence $3K_1$ is not an alternative, we have to substitute
 the other regular graph, $K_3$, which has internal degree 2.
In this way we obtain the degrees
 $$
   11, \ 12, \ 14
 $$
 which cannot be coupled with the case $(12, 13, 14)$ of $16=6+6+4$.

In $15=6+5+4$ the `5' class means $C_5$ with internal degree 2 and external degree 10, i.e.\ degree 12 in $F$.
Therefore the `4' class with external degree 11 must be $C_4$ with internal degree 2 and total degree 13.
The external degree for `6' is 9, hence internal degree 3 is infeasible, thus we have to substitute $2K_3$
 which leads to degree 11 and in this way we obtain the degree set
 $$
   11, \ 12, \ 13.
 $$
From this, it is clear that $16+15$ cannot occur, and even $15+15$ would be impossible.
 (In fact, degree 12 appears in all the three types above, and any two types have
   two values in common.)
\qed

\brm
 The construction on 30 vertices is another example of the mixed principle
  as described in Section \ref{ss:mixedconst}.
 Here we start from the graph $H=2K_3+e$, two vertex-disjoint triangles
  connected by just one edge $e$.
 Although this $H$ is not paw-free, still does not contain a singular paw;
  and its complement $\overline{H}\cong K_{3,3}-e$ is paw-free.
 Then the two ends of the edge $e$ can be viewed together
  as one partition class, while the other classes are singletons.
 Each end of $e$ has two neighbors in $H$ and this yields two neighbor classes for
  the corresponding subsets after substitution.
 In case of the paw, two classes of order 6 with identical neighborhood
  may occur because their internal degree can (and should) be distinct.
\erm

\subsection{The 4-cycle $C_4$}

In case of $C_4$, which seems most problematic among the small graphs,
 we can derive lower and upper bounds which are quite close to each other,
 but still the exact value of $\Rs(C_4)$ is unknown.

Note that the 4-cycle has $\R{C_4}=6$, and its two R-graphs of order 5
 are $C_5$ and the bull.
  Neither of them is Ramsey-stable.
Indeed, removing a vertex from $C_5$ we obtain $P_4$, which is extendable
 not only to $C_5$ itself, but also to the bull.
Similarly, removing the degree-2 vertex from the bull we obtain $P_4$ which is
 extendable to $C_5$.
Moreover, the removal of a pendant vertex from the bull yields the paw, which can be
 extended to the bull in two different ways.
Also, removing a vertex of degree 3 we obtain $P_3\cup K_1$, whose extension to the
 bull fixes an edge to the isolated vertex, and another edge to the middle of $P_3$,
 but the last edge can go to either end of $P_3$.

\bp
 $24 \le \Rs(C_4) \le 26$.
\ep

 \pf
The upper bound is a consequence of Corollary \ref{c:square}.
 For the lower bound we construct an SR-graph on 23 vertices.
Let us take the bull as host graph $H$, labeling
 its vertices as $v_1,\dots,v_5$ where $\{v_1,v_2,v_3\}$ induces a triangle, and the
 two pendant edges are $v_1v_4$ and $v_3v_5$.
Let us substitute graphs $F_i$ for $v_i$ such that $F_1\cong K_3$
 and $F_i\cong C_5$ for all $2\le i\le 5$.
All internal degrees are equal to 2, and the external degrees are 15 in $F_1$,
 8 in $F_2$, 13 in $F_3$, 3 in $F_4$, and 5 in $F_5$.
Neither $F$ nor its complement contains any singular copy of $C_4$, hence $\Rs(C_4)\ge 24$.
 \qed

\subsection{Small graphs with isolates}

In this last of the subsections devoted to small graphs we give the values
 for those graphs of order four which have isolated vertices.
There are three such graphs: $K_2\cup 2K_1$, $P_3\cup K_1$, and $K_3\cup K_1$.
 Note that some lower bounds can easily be obtained from above:
  \begin{itemize}
    \item Theorem \ref{th:2} (with reference also to Remark \ref{rm:mon})
     implies
      $$
        \Rs(P_3\cup K_1) \ge \Rs(K_2\cup 2K_1) \ge 10 .
      $$
    \item The construction of Theorem \ref{th:7} yields
      $$
        \Rs(K_3\cup K_1) \ge 22 .
      $$
  \end{itemize}
We prove that these bounds are tight.

\bp
 $\Rs(P_3\cup K_1) = \Rs(K_2\cup 2K_1) = 10$.
\ep

 \pf
In every graph $G$ with 10 vertices there exists a singular subgraph
 of order four.
It necessarily contains a $P_3$ or its complement, which can be extended
 to a singular $P_3\cup K_1$.
Thus $\Rs(P_3\cup K_1) \le 10$.
 \qed

\thm
 $\Rs(K_3\cup K_1) = 22$.
\ethm

 \pf
Suppose for a contradiction that $F$ is an SR-graph of order $n\geq 22$
  for $K_3\cup K_1$.
We know that a singular $K_3$ occurs in $F$ (or in its complement),
 say it has the vertices $v_1,v_2,v_3$.
If the degrees of this $K_3$ are all distinct, say $d_1<d_2<d_3$, then it would be
 extendable to a singular $K_3\cup K_1$ unless all vertices of $F$ have their
 degree from $\{d_1,d_2,d_3\}$.
But then a degree class would have at least eight vertices, so that $F$ would contain
 even a singular $K_3\cup 5K_1$.

Hence suppose that the three vertices of any singular $K_3$ have the same
 degree in $F$. 
Since $\R{K_3}=6$, there can be at most five degree classes, and we easily find
 a singular $K_3\cup K_1$ unless all degree classes have at most five vertices
 and the degree class(es) inducing a triangle have exactly three vertices.
In particular, $\{v_1,v_2,v_3\}$ itself is a degree class, moreover its
 complementary 19 or more vertices form only four degree classes.
Now we can only have the following possibilities:
 \begin{itemize}
   \item $n=22=3+4+5+5+5$,
   \item $n=23=3+5+5+5+5$.
 \end{itemize}
And then the degree counting method in the proof of Theorem \ref{th:7} can be
 repeated for these two cases without any changes, leading to the contradiction
  $\delta(F)\geq 9$ for $n=22$ and $\delta(F)\geq 10$ for $n=23$.
 \qed

\section{Stars of any size}
   \label{s:star}

The star with $s$ edges, $K_{1,s}$, is an easy case concerning Ramsey numbers;
  cf.\ e.g.\ Section 5.5 of \cite{R-surv}:
 \begin{itemize}
   \item if $s$ is odd, then $\R{K_{1,s}}=2s$, and the extremal R-graphs are
     precisely the $(s-1)$-regular graphs of order $2s-1$;
   \item if $s$ is even, then $\R{K_{1,s}}=2s-1$, and the extremal R-graphs are
      the graphs of order $2s-2$ with minimum degree at least $s-2$ and
     maximum degree at most $s-1$.
      In particular, the largest \textit{regular} R-graphs are those graphs
     of order $2s-2$ which are $(s-2)$-regular or $(s-1)$-regular.
 \end{itemize}

It turns out that the parity of $s$ is essential also with respect to $\Rs$.
The case of even $s$ is simpler, the quadratic upper bound always is tight.

\bp
 If $s$ is even, then
  $$
    \Rs(K_{1,s}) = (\R{K_{1,s}} - 1)^2 + 1 = (2s-2)^2 + 1 .
  $$
\ep

\pf
 The upper bound $(2s-2)^2 + 1$ follows from Corollary \ref{c:square}.
For the same lower bound we construct an RS-graph $F$ of order $(2s-2)^2$
 on vertex set
  $$
    V = (A_1 \cup \cdots \cup A_{s-1}) \cup (B_1 \cup \cdots \cup B_{s-1})
  $$
 where all sets $A_i$ and $B_i$ are mutually disjoint, each of cardinality $2s-2$.
The edges of $F$ are defined as follows:
 \begin{itemize}
   \item each $A_i$ induces an $(s-2)$-regular graph;
   \item each $B_i$ induces an $(s-1)$-regular graph;
   \item there is no edge between $A_i$ and $A_j$ for $i\ne j$;
   \item there is no edge between $B_i$ and $B_j$ for $i\ne j$;
   \item every $A_i$ and $B_j$ are completely adjacent for $i\ne j$;
   \item the sets $A_i$ and $B_i$ are adjacent by a $(2i)$-regular bipartite
     graph, for all $1\le i\le s-1$.
 \end{itemize}
Then, both in $F$ and in $\overline{F}$, each vertex $v$ is adjacent to at most
 $s-1$ vertices of distinct degrees different also from the degree of $v$;
 and to at most $s-1$ vertices whose degree is equal to that of $v$.
Thus, $F$ is an RS-graph for $K_{1,s}$.
 \qed

\bsk

 The case of odd $s$ is more complicated.
The quadratic upper bound is never attained, although the singular Ramsey number
 is not far from it.
For the tightness of the lower bound we give two very different constructions,
 with the purpose to indicate that --- contrary to $\R{K_{1,s}}$ --- the extremal
  graphs for $\Rs(K_{1,s})$ may be quite hard to characterize.

\thm
 If $s$ is odd, then
  $$
    \Rs(K_{1,s}) = (\R{K_{1,s}} - 1)^2 + 1 - (2s-2) = (2s-1)(2s-2) + 2 .
  $$
\ethm

 \nin
{\bf Proof of the Upper Bound.}
Suppose for a contradiction that there exists an SR-graph $F$ of order
 at least $(2s-1)(2s-2) + 2$ for $K_{1,s}$.
  Denote its degree classes by $V_1,\dots,V_m$.
We know that $m\le 2s-1$, and also $|V_i|\le 2s-1$ for all $1\le i\le m$.
 Hence $m = 2s-1$ must hold.
Since all R-graphs of order $2s-1$ are regular, all of them are Ramsey-stable, due to
 Remark \ref{r:regRstab}.
Thus, by the Regular Substitution Lemma, each $V_i$ induces a regular R-graph,
 and two distinct $V_i,V_j$ are either completely adjacent or completely
 nonadjacent,
From this we obtain that the structure of adjacencies between the degree classes
 is an $(s-1)$-regular graph of order $2s-1$, therefore
 \begin{itemize}
   \item every external degree is at most $(s-1)(2s-1)$, and every internal
    degree is at most $s-1$, therefore the maximum degree of $F$ is
    at most $2s(s-1)$ and the minimum degree
    cannot be larger than $2s(s-1) - (2s-2) = 2(s-1)^2$.
 \end{itemize}
Let us define $r_i := (2s-1) - |V_i|$ for $i=1,\dots, m$.
 Then the internal degree inside $V_i$ is at least $(s-1) - r_i$.
Moreover, if a $V_j$ is adjacent to $V_i$, then it contributes to the
 external degree of every $v\in V_i$ with exactly $(2s-1) - r_j$.
It follows that
 \begin{itemize}
   \item the minimum degree is at least $2s(s-1) - \sum_{i=1}^{2s-1} r_i$
 \end{itemize}
  from where we obtain that
 $$
   \sum_{i=1}^{2s-1} r_i \ge 2s-2
 $$
  and
 $$
   |V(F)| = (2s-1)^2 - \sum_{i=1}^{2s-1} r_i \le (2s-1)^2 - (2s-2) < (2s-1)(2s-2) + 2 ,
 $$
  a contradiction.

\bsk

 \nin
{\bf Proof of the Lower Bound.}
 Let $s=2t+1$; then $\R{K_{1,s}} - 1 = 2s-1 = 4t+1$.
We start with the graph $H = (C_{4t+1})^t$ which has vertices $x_0,x_1,\dots,x_{4t}$
 and edges $x_ix_{i\pm j}$ for $j=1,\dots,t$, where subscript addition is taken
 modulo $4t+1$.
For each $x_i$ we substitute a $d_i$-regular graph $G_i$ with vertex set $V_i$
  in the following way:
 \begin{itemize}
   \item $2t+1$ sets $|V_1|=\dots=|V_{2t+1}|=4t+1$ and $d_1=\dots=d_{2t+1}=2t$;
   \item $t$ sets $|V_{2t+2}|=\dots=|V_{3t+1}|=4t$ and $d_{2t+2}=\dots=d_{3t+1}=2t-1$;
   \item $t-1$ sets $|V_{3t+2}|=\dots=|V_{4t}|=4t$ and $d_{3t+2}=\dots=d_{2t+1}=2t$;
   \item $1$ set $|V_{0}|=2t$ and $d_{0}=t$.
 \end{itemize}
Recall that $V_i$ is adjacent to $V_{i-t}, V_{i-t+1}, \dots, V_{i-1}, V_{i+1},
 V_{i+2}, \dots,V_{i+t}$, where subscript addition is taken modulo $4t+1$.
Then the obtained degrees --- more precisely their differences from the maximum of the
 internal / external / total degree --- can be summarized as shown in Table \ref{tab1}.

\renewcommand{\arraystretch}{1.5}
\begin{table}[htp]
	\begin{center}%
	{
	 \tiny
\begin{tabular}{|c||c|c|c|c|c|c||c|c|c|c|c|c|c||c|}
\hline 
  & $V_{1}$ & $\dots$ & $V_{t}$ & $V_{t+1}$ & $\dots$ & $V_{2t+1}$ & $V_{2t+2}$ & $\dots$ & $V_{3t}$ & $V_{3t+1}$ & $V_{3t+2}$ & $\dots$ & $V_{4t}$ & $V_{0}$ \\ 
\hline 
\hline 
SZ & 0 & $\dots$ & 0 & 0 & $\dots$ & 0 & 1 & $\dots$ & 1 & 1 & 1 & $\dots$ & 1 & $2t+1$ \\ 
\hline 
ID & 0 & $\dots$ & 0 & 0 & $\dots$ & 0 & 1 & $\dots$ & 1 & 1 & 0 & $\dots$ & 0 & $t$ \\ 
\hline 
ED & $3t$ & $\searrow$ & $2t+1$ & 0 & $\nearrow$ & $t$ & $t$ & $\nearrow$ & $2t-2$ & $4t-1$ & $4t-1$ & $\searrow$ & $3t+1$ & $t$ \\ 
\hline 
TD & $3t$ & $\searrow$ & $2t+1$ & 0 & $\nearrow$ & $t$ & $t+1$ & $\nearrow$ & $2t-1$ & $4t$ & $4t-1$ & $\searrow$ & $3t+1$ & $2t$ \\ 
\hline 
\end{tabular} 
  }
\end{center}
\caption{Some parameters of the extremal construction for unrestricted $s=2t+1$.
  SZ = \{$(2s-1)$ minus size\} = $4t+1-|V_i|$\,; \ ID = \{$(s-1)$ minus internal degree\} = $2t-d_i$\,; \
  ED = \{$(s-1)(2s-1)$ minus external degree\}\,; \ TD= \{$2s(s-1)$ minus total degree\}\,; \
  {\footnotesize $\searrow \;, \; \nearrow$} = decreasing / increasing by 1 in each step. }\label{tab1}
\end{table}
Then the degrees range between $2s(s-1) - 4t = 2s(s-1) - (2s-2)$ and $2s(s-1)$,
 and the number of vertices is $16t^2 + 4t + 1 = (2s-1)(2s-2) + 1$.
Both in the graph and in its complement, each vertex has neighbors only in $s-1$ other degree classes,
 and at most $s-1$ neighbors in its degree class.
Hence we have an SR-graph of the required order.

\bsk

 \nin
  {\bf Alternative construction for \boldmath${s=4q+1}$.}
The basic structure $H = (C_{4t+1})^t =  (C_{8q+1})^{2q}$ remains the same, but the
 size distribution of substituted R-graphs will be substantially different: they will have
 almost equal sizes, rather than involving a very small degree class.
We need a construction on
 $(2s-1)^2 - (2s-2) = (8q+1)^2 - 8q$ vertices.
This will be achieved by taking $4q+1$ degree classes of size $8q+1$, moreover
 $2q$ classes of size $8q$, and $2q$ classes of size $8q-2$.

We use the symbol $G_{p}^{d}$ to denote any $d$-regular graph on $p$ vertices.
 Such graphs exist whenever $p>d\ge 0$ and $pd$ is even.
In the construction below, the actual structure of a $G_{p}^{d}$ will be irrelevant, one may take
 different graphs for different appearances of the same pair $(p,d)$.
Using the notation $V_i$ and $G_i$ in the sense as above, we now define:

 \begin{itemize}
   \item for every $i$ in the range $2q\le i\le 6q$ we take
    $|V_{i}|=8q+1$, and let each $G_i$ be a $G_{2s-1}^{s-1}$;
   \item with the only one exception of $V_{8q}$,
    for all $1\le i\le q$ we take $|V_{4q\pm(2q+2i)}|=8q$, and let each $G_i$
    be a $G_{2s-2}^{s-2}$;
   \item with the only one exception of $V_{2q-1}$,
    for all $1\le i\le q$ we take $|V_{4q\pm(2q+2i-1)}|=8q-2$, and let each $G_i$
    be a $G_{2s-4}^{s-2}$;
   \item for the two exceptional cases we take $|V_{8q}|=8q-2$ with
    $G_{8q}=G_{2s-4}^{s-3}$ and $|V_{2q-1}|=8q$ with $G_{2q-1}=G_{2s-2}^{s-1}$.
 \end{itemize}
The maximum degree occurs at the vertices of $V_{4q}$\,: they have
 internal degree $s-1$, external degree $(s-1)(2s-1)$, and total degree $2s(s-1)$.
Relevant parameters of vertices in the other degree classes are
 summarized in Table~\ref{tab2}.
One can check that each $V_i$ has a distinct degree, and consequently we obtained
 an SR-graph of maximum order.
  \qed

\noindent
\renewcommand{\arraystretch}{1.5}
\begin{table}[htp]
	\begin{center}%
	{\scriptsize
\begin{tabular}{|c||c|c|c|c|c|c|c|c|c|c|}
\hline 
  & $V_{4q}$ & $V_{4q-1}$ & $V_{4q-2}$ & $V_{4q-3}$ & $V_{4q-4}$ & $\dots$ & $V_{2q+3}$ & $V_{2q+2}$ & $V_{2q+1}$ & $V_{2q}$ \\ 
\hline 
\hline 
SZ & 0 & 0 & 0 & 0 & 0 & $\dots$ & 0 & 0 & 0 & 0 \\ 
\hline 
ID & 0 & 0 & 0 & 0 & 0 & $\dots$ & 0 & 0 & 0 & 0 \\ 
\hline 
ED & 0 & 1 & 2 & 5 & 6 & $\dots$ & $4q-7$ & $4q-6$ & $4q-3$ & $4q-2$ \\ 
\hline 
TD & 0 & 1 & 2 & 5 & 6 & $\dots$ & $4q-7$ & $4q-6$ & $4q-3$ & $4q-2$ \\ 
\hline 
\end{tabular} 
 \\ ~~ \medskip ~~ \\
\begin{tabular}{|c||c|c|c|c|c|c|c|c|c|}
\hline 
  & $V_{4q+1}$ & $V_{4q+2}$ & $V_{4q+3}$ & $V_{4q+4}$ & $\dots$ & $V_{6q-3}$ & $V_{6q-2}$ & $V_{6q-1}$ & $V_{6q}$ \\ 
\hline 
\hline 
SZ & 0 & 0 & 0 & 0 & $\dots$ & 0 & 0 & 0 & 0 \\ 
\hline 
ID & 0 & 0 & 0 & 0 & $\dots$ & 0 & 0 & 0 & 0 \\ 
\hline 
ED & 3 & 4 & 7 & 8 & $\dots$ & $4q-5$ & $4q-4$ & $4q-1$ & $4q+2$ \\ 
\hline 
TD & 3 & 4 & 7 & 8 & $\dots$ & $4q-5$ & $4q-4$ & $4q-1$ & $4q+2$ \\ 
\hline 
\end{tabular} 
 \\ ~~ \medskip ~~ \\
\begin{tabular}{|c||c|c|c|c|c|c|c|c|c|c|c|}
\hline 
  & $V_{2q-1}$ & $V_{2q-2}$ & $V_{2q-3}$ & $V_{2q-4}$ & $V_{2q-5}$ & $V_{2q-6}$ & $\dots$ & $V_{3}$ & $V_{2}$ & $V_{1}$ & $V_{0}$ \\ 
\hline 
\hline 
SZ & 1 & 1 & 3 & 1 & 3 & 1 & $\dots$ & 3 & 1 & 3 & 1 \\ 
\hline 
ID & 0 & 1 & 1 & 1 & 1 & 1 & $\dots$ & 1 & 1 & 1 & 1 \\ 
\hline 
ED & 4q & 4q+3 & 4q+2 & 4q+7 & 4q+6 & 4q+11 & $\dots$ & $8q-10$ & $8q-5$ & $8q-6$ & $8q-1$ \\ 
\hline 
TD & 4q & 4q+4 & 4q+3 & 4q+8 & 4q+7 & 4q+12 & $\dots$ & $8q-9$ & $8q-4$ & $8q-5$ & $8q$ \\ 
\hline 
\end{tabular} 
 \\ ~~ \medskip ~~ \\
\begin{tabular}{|c||c|c|c|c|c|c|c|c|c|c|c|}
\hline 
  & $V_{6q+1}$ & $V_{6q+2}$ & $V_{6q+3}$ & $V_{6q+4}$ & $V_{6q+5}$ & $V_{6q+6}$ & $\dots$ & $V_{8q-3}$ & $V_{8q-2}$ & $V_{8q-1}$ & $V_{8q}$ \\ 
\hline 
\hline 
SZ & 3 & 1 & 3 & 1 & 3 & 1 & $\dots$ & 3 & 1 & 3 & 3 \\ 
\hline 
ID & 1 & 1 & 1 & 1 & 1 & 1 & $\dots$ & 1 & 1 & 1 & 2 \\ 
\hline 
ED & 4q & 4q+5 & 4q+4 & 4q+9 & 4q+8 & 4q+13 & $\dots$ & $8q-8$ & $8q-3$ & $8q-4$ & $8q-3$ \\ 
\hline 
TD & 4q+1 & 4q+6 & 4q+5 & 4q+10 & 4q+9 & 4q+14 & $\dots$ & $8q-7$ & $8q-2$ & $8q-3$ & $8q-1$ \\ 
\hline 
\end{tabular} 
  }
\end{center}
\caption{Parameters in the following ranges: $4q-1\ge i\ge 2q$\,; \,$4q+1\le i\le 6q$\,;
 \,$2q-1\ge i\ge 0$\,; \,$6q+1\le i\le 8q$. Notations SZ, ID, ED, TD are the same as in
  Table \ref{tab1}. }\label{tab2}
\end{table}

\newpage

\section{Asymptotics for singular Tur\'an numbers}
   \label{s:tur}

In this section we present estimates on the singular Tur\'an numbers $\Ts(n,F)$,
  and compare them to the classical Tur\'an number $\mbox{\rm ex}(n,K_s)$, which is
  the maximum number of edges in a graph of order $n$ not containing
  a complete subgraph of order $s$.
 Assume\footnote{If $p=2$, then either $F=2K_1$ which is a singular subgraph
  of every graph with at least two vertices --- hence $\Ts(n,F)$ is
  meaningless --- or $F=K_2$ and $\Ts(n,F)=0$ for all~$n$.
 The situation is similar if $\chi(F)=1$, i.e.\ $F=pK_1$, in which case
  $\Ts(n,F)$ cannot be defined for $n>(p-1)^2$.}
  that $F$ has order $p:=|V(F)|\ge 3$ and chromatic number $q:=\chi(F)\ge 2$.

We begin with two general constructions, providing lower bounds on $\Ts(n,F)$.

Let us assume that $n$ is a multiple of $q-1$; regarding lower bounds
 for other orders we refer to the simple fact that
  $$\Ts(n,F) \ge \Ts\!\left((q-1)\left\lfloor \frac{n}{q-1} \right\rfloor,F\right) .$$
This will give a fairly good approximation because the difference between the
 numbers of edges for two consecutive multiples of $q-1$ will be $O(n)$ only,
 while $\Ts$ will be shown to grow with a quadratic function of $n$.
Even in a more general setting where $F$ is not a fixed graph and $q$ varies, say
 $F=K_{\lfloor \sqrt{n} \rfloor}$,
 the difference between the numbers of edges will grow with $O(qn)=o(n^2)$,
 which is negligible compared to $\Ts(n,F)$.

The higher structure of both constructions is a partition of the $n$-element vertex set
 into $p-1$ classes $V_1,V_2, \dots , V_{p-1}$, where
 \begin{itemize}
   \item each $|V_i|$ is a multiple of $q-1$, and
   \item each $V_i$ induces a Tur\'an graph for $K_q$, i.e., the subgraph induced by $V_i$
    in the graph of order $n$ under construction is a complete multipartite graph
    with $q-1$ vertex classes $U_{i,1},U_{i,2}, \dots , U_{i,q-1}$ of equal size,
  $$|U_{i,1}| = |U_{i,2}| = \ldots = |U_{i,q-1}| ,$$
     each $U_{i,j}$ is an independent set, and any two of them are
     completely adjacent.
 \end{itemize}
In both constructions the degree classes will be $V_1,V_2, \dots , V_{p-1}$.

 \bcs   \label{cs:1}
Choose the sizes of the degree classes $V_i$ in such a way that
  $$|V_{1}| < |V_{2}| < \cdots < |V_{p-1}|$$
 holds, and under this condition $|V_{1}|$ is as large as possible,
  whereas $|V_{p-1}|$ is as small as possible.
 \ecs

Since the sequence $|V_{1}|, |V_{2}|, \dots , |V_{q-1}|$ is
 strictly increasing, we must have
 $|U_{i,j}| \geq |U_{1,j}|+i-1$ for all $1\le i<p-1$
 (and all $1\le j\le q-1$).
Then the requirement on $|V_{1}|$ and $|V_{p-1}|$ means that
 we need to maximize $|U_{1,1}|$ subject to
 $$
   (p-1)(q-1)|U_{1,1}| + (q-1) \sum_{i=2}^{p-1} (i-1) ~ \le ~ n ,
 $$
  from where we obtain that
   $|U_{1,j}| \approx \frac{n}{(p-1)(q-1)} - \frac{p}{2}$
    and $|U_{p-1,j}| \approx \frac{n}{(p-1)(q-1)} + \frac{p}{2}$.
In fact either $|U_{p-1,1}|=|U_{1,1}|+p-2$ or
 $|U_{p-1,1}|=|U_{1,1}|+p-1$.
By construction we also have:
 \begin{itemize}
   \item the vertices of $V_i$ have degree $n - |U_{i,1}|$.
 \end{itemize}
This implies that the degree sets are indeed the classes $V_i$, and
 two types of singular subgraphs can occur:
 \begin{itemize}
   \item subgraphs of a $V_i$, thus having chromatic number less than $q$;
   \item subgraphs with at most one vertex in each $V_i$, thus having order
    less than $p$.
 \end{itemize}
It follows that the constructed graph does not contain any singular subgraph
 isomorphic to $F$.

Let us compare the number of edges with that in the Tur\'an graph for
 $K_{(p-1)(q-1)+1}$.

\bp   \label{p:turlo}
 Let $F$ be a graph with $p\ge 3$ vertices and chromatic number $q\ge 2$.
 If $n$ is a multiple of $q-1$, then
  $$
    \mbox{\rm ex}(n,K_{(p-1)(q-1)+1}) - \Ts(n,F) \le cqp^3 .
  $$
   for a constant $c$.
 If $n\equiv r ~ (\mbox{\rm mod } (q-1))$ with $r\ne 0$, then
  $$
    \mbox{\rm ex}(n,K_{(p-1)(q-1)+1}) - \Ts(n,F) \le O(rn) .
  $$
\ep

 \pf
 Suppose first that $n$ is divisible by $q-1$.
From the graph obtained in Construction \ref{cs:1} we can obtain the Tur\'an graph
 if, for every $1\le i\le\frac{p-1}{2}$ and every $1\le j\le q-1$, we replace the
 vertex classes $U_{i,j}$ and $U_{p-i,j}$ with two classes (independent sets) of sizes
 $\left\lfloor \frac{|U_{i,j}|+|U_{p-i,j}|}{2} \right\rfloor$
 and $\left\lceil \frac{|U_{i,j}|+|U_{p-i,j}|}{2} \right\rceil$.
Due to the identity $(x-a)(x+a)=x^2-a^2$, this operation increases the number of
 edges proportionally to $(p/2-i)^2$,
  because the subgraph induced by $U_{i,j} \cup U_{p-i,j}$ remains
 a complete bipartite graph on exactly the same set of vertices and
  with an unchanged number of edges to its exterior.
There are $q-1$ choices for $j$, and $i$ runs from 1 to $\lfloor (p-1)/2 \rfloor$,
 hence the total difference grows with the order of $qp^3$.

If $n=t(q-1)+r$ with $r\ne 0$, then we supplement the construction with $r$
 isolated vertices, hence no singular $F$ will arise while the number of edges
 does not decrease (actually remains unchanged).
On the other hand, the Tur\'an function clearly satisfies the inequality
 $\mbox{\rm ex}(n,H) - \mbox{\rm ex}(n-r,H) < rn$ for every graph $H$ and
  all natural numbers $n$ and $r$.
 \qed

 \bcs
For the sake of simpler description assume that $n$ is a multiple of $(p-1)(q-1)$,
 with $q\ge 2$ and $p\ge 4$.
We define all $U_{i,j}$ to have the same size ($1\le i\le p-1$, $1\le j\le q-1$),
 i.e.\ $|U_{i,j}| = \frac{n}{(p-1)(q-1)}$, each of them being an independent set;
 hence in particular $|V_{i}| = \frac{n}{p-1}$,
  where $V_i=\bigcup_{j=1}^{q-1} U_{i,j}$.
  Start with complete bipartite graphs between any two $U_{i_1,j_1},U_{i_2,j_2}$.
Represent the sets $V_{i}$ with single vertices $v_{i}$, and view them as
 the vertices of $K_{p-1}$.
It was proved by Chartrand et al.\ \cite{Ch+} that the edges of $K_{p-1}$
 can be assigned with integer weights from $\{1,2,3\}$ in such a way that the
 weighted degrees of the vertices become mutually distinct.
Now, for each vertex pair,
 \begin{itemize}
   \item if the weight of $v_{i_1}v_{i_2}$ is 1, keep complete
    adjacency between $V_{i_1}$ and $V_{i_2}$;
   \item if the weight of $v_{i_1}v_{i_2}$ is 2, remove a
    perfect matching between $V_{i_1}$ and $V_{i_2}$;
   \item if the weight of $v_{i_1}v_{i_2}$ is 3, remove a
    2-factor between $V_{i_1}$ and $V_{i_2}$.
 \end{itemize}
 \ecs

By construction, the degree classes are the sets $V_i$,
 hence the graph does not contain any singular subgraph with $p$ vertices
 and chromatic number $q$; and the number of removed edges is at most
 $(p-2) n$.
In fact, applying the results of \cite{RSA95},
 this upper bound can be reduced to
 $(p/2 + c) n$, where $c$ is a universal constant for all $n$ and $p$,
 which is tight apart from the actual value of $c$.

Next, we prove an upper bound which shows that the constructions above give
 tight asymptotics on $\Ts(n,F)$ for every fixed graph $F$ as $n$ gets large.

 \thm   \label{th:asympt}
If $F$ is a graph with $p\ge 3$ vertices and chromatic number $q\ge 2$, then
  $$
    \Ts(n,F) \le \mbox{\rm ex}(n,K_{(p-1)(q-1)+1}) + o(n^2) .
  $$
Moreover, for the complete graph $K_p$ (i.e., $q=p$) we have
  $$
    \Ts(n,K_p) \le \mbox{\rm ex}(n,K_{(p-1)^2+1}) .
  $$
Both upper bounds are asymptotically sharp as $n\to\infty$.
 \ethm

\pf
 We begin with the inequality for $K_p$, as it is much simpler to prove.
If a graph $G$ of order $n$ has more than $\mbox{\rm ex}(n,K_{(p-1)^2+1})$
 edges, then by definition it contains a complete subgraph on
 $(p-1)^2+1$ vertices.
Among them, $p$ have the same degree in $G$ or $p$ have mutually
 distinct degrees.
Thus, $G$ contains $K_p$ as a singular subgraph, which implies that
 $\Ts(n,K_p)$ cannot be that large.

In the general case let us asssume that $G$ is a graph of order $n$,
 having as many as $\mbox{\rm ex}(n,K_{(p-1)(q-1)+1}) + \epsilon n^2$ edges.
We are going to apply the Erd\H os--Stone theorem \cite{ES46},
 which states that for any fixed $\epsilon>0$ a graph with $n$ vertices
  and $\mbox{\rm ex}(n,K_{s}) + \epsilon n^2$ edges contains not only
  a $K_s$ but also a complete multipartite graph with $s$ vertex classes
  with $t$ vertices in each class;
 here $t$ can be taken any large as $n$ increases.\footnote{For our
  purpose with a fixed $F$, the classical theorem by Erd\H os and Stone
  from 1946 is sufficiently strong.
 An improved numerical estimate on $t$ was derived three decades later
  by Bollob\'as et al.\ in \cite{BES76}, and finally Chv\'atal and
  Szemer\'edi proved in \cite{ChSz81} that $t$ grows as fast as
  $c\frac{\log n}{\log 1/\epsilon}$.
 This version is useful when one takes a sequence of graphs $F$ whose
  orders tend to infinity as $n$ gets large but does not exceed
  $c'\sqrt{\log n}$ for a small constant $c'$, e.g.\ in case of
  $F=K_{\lceil\log\log n\rceil}$.}
We take $s=(p-1)(q-1)+1$ and assume that $n$ is large enough
 to ensure that also $t$ is sufficiently large, say $t\ge p^2$.

Let $A_1,\dots,A_{(p-1)(q-1)+1}$ be the vertex classes of a
 $K_{t,\dots,t}\subset G$.
Each $A_i$ contains a singular $B_i\subset A_i$ with $|B_i|\ge\sqrt{t}\ge p$,
 with vertices whose degrees are all equal or all distinct in $G$.

If the degrees are all distinct in at least $p$ of the sets $B_i$,
 then we can sequentially select one vertex from each $B_i$ such that
 in each step the degree of the selected vertex is distinct from
 all previously selected ones.
This yields a singular $K_p\subset G$, thus also $F$ occurs as a
 singular subgraph of $G$.

Suppose that there are only $h$ sets $B_i$ (where $0\le h\le p-1$)
 inside which the degrees are distinct.
We assume that these classes are the ones with largest subscripts,
 namely $B_{(p-1)(q-1)-h+2},\dots,B_{(p-1)(q-1)+1}$.
Then we obtain a sequence $d_1,d_2,\dots,d_{(p-1)(q-1)-h+1}$
 where $d_i$ is the degree of all vertices in $B_i$.
If this sequence contains $q$ equal terms, then from the corresponding
 sets we can select the vertices for the color classes of $F$ in a
 proper $q$-coloring, thus $F$ is a singular subgraph of $G$ with
 all degrees equal.
Else every value occurs at most $q-1$ times, hence the sequence
 contains at least
  $z:=\left\lceil \frac{(p-1)(q-1)-h+1}{q-1} \right\rceil$
 mutually distinct terms, which can be supplemented at least with further
 $\min(h,p-z)$ distinct degrees from the last $h$ sets $B_i$.
Now we have
  $$
    \left\lceil \frac{(p-1)(q-1)-h+1}{q-1} \right\rceil + h
      = p-2 + h + \left\lceil \frac{q-h}{q-1} \right\rceil
        \ge p ,
  $$
 with equality only if $h=0$ or $h=1$ or $q=h=2$.
Consequently, $K_p$ occurs as a singular subgraph of $G$ with
 all degrees distinct.

\msk

Asymptotic tightness follows from the constructions described above,
 for both cases.
\qed

\bsk

For the case $p = q = 3$ and $n\equiv 2 \ (\mbox{\rm mod }4)$
 we obtained an exact result.

\bc
 If $F=K_3$ (i.e., $p=q=3$) and $n\equiv 2 \ (\mbox{\rm mod }4)$, then
$$
    \Ts(n,K_3) = \mbox{\rm ex}(n,K_5)
      = \frac{3}{8} \, n^2 - \frac{1}{2} .
$$
\ec

\pf
 Assume that $n=4h+2$.
Then the Tur\'an graph for $K_5$ is the complete 4-partite graph
 in which the vertex classes have respective cardinalities
 $h,h,h+1,h+1$.
The first $2h$ vertices have degree $3h+2$, while the last $2h+2$
 vertices have degree $3h+1$.
Hence there are only two degree classes, each of them inducing
 a complete bipartite graph, therefore the graph certainly is $K_3$-free.
Thus no singular $K_3$ occurs, implying $\Ts(n,K_3) \ge \mbox{\rm ex}(n,K_5)$.
Also the reverse inequality is valid, by Theorem \ref{th:asympt}.
 \qed

\bsk

We close this section with some fairly tight estimates for $K_3$.

\bp   \label{p:trian}
 For $F=K_3$ and $n\ge 3$ we have the following inequalities.
  \begin{itemize}
    \item[$(i)$] If $n\equiv 0 \ (\mbox{\rm mod }4)$, then
      $
        \frac{3}{8} \, n^2 - 2
         \le \Ts(n,K_3) \le \frac{3}{8} \, n^2 - 1 .
      $
    \item[$(ii)$] If $n\equiv 1 \ (\mbox{\rm mod }4)$, then
      $
        \frac{3}{8} \, n^2 - \frac{1}{4} \, n - \frac{1}{8}
         \le \Ts(n,K_3) \le \frac{3}{8} \, n^2 - \frac{11}{8} .
      $
    \item[$(iii)$] If $n\equiv 3 \ (\mbox{\rm mod }4)$, then
      $
        \frac{3}{8} \, n^2 - \frac{1}{4} \, n - \frac{13}{8}
         \le \Ts(n,K_3) \le \frac{3}{8} \, n^2 - \frac{11}{8} .
      $
  \end{itemize}
\ep

 \pf
In all cases, the claimed upper bound is a Tur\'an number
 minus 1, namely $\mbox{\rm ex}(n,K_5) - 1$.
Its validity follows from Theorem \ref{th:asympt} by the further
 observation that the corresponding Tur\'an graphs are unique
 and each of them contains a singular $K_3$.
For the lower bounds we give constructions as follows.

\msk

$(i)$\quad 
 This is a particular case of Construction \ref{cs:1}, with
  $|U_{1,1}| = |U_{1,2}| = \frac{1}{4} \, n - 1$ and
  $|U_{2,1}| = |U_{2,2}| = \frac{1}{4} \, n + 1$.

\msk

$(ii)$\quad 
 Start with the complete 4-partite graph with equal vertex classes
  of size $\frac{1}{4} \, (n-1)$, and join a new vertex, say $z$,
  to all vertices of two classes.
 Denoting $n=4h+1$, the two classes adjacent to $z$ have degree $3h+1$,
  the other two classes have degree $3h$, and $z$ has degree $2h$.
 There is no singular $K_3$ because there are only three distinct
  degrees (degree $2h$ occurring only on $z$),
   each degree class is triangle-free, and in every triangle
   containing $z$ the other two vertices have degree $3h+1$.

\msk

$(iii)$\quad 
 Assume that $n=4h+3$.
 Start with the optimal construction for $n-1$, that is the
  complete 4-partite graph with vertex classes of respective
   sizes $h,h,h+1,h+1$.
 Similarly to the case of $(ii)$ join a new vertex $z$ to the
  $2h$ vertices of the two smaller classes.
 Then $2h$ vertices have degree $3h+3$, $2h+2$ vertices have degree $3h+1$,
  and $z$ has degree $2h$ alone, with all its neighbors having
  degree $3h+3$.
\qed

\bsk

The principle of these constructions can also be applied to obtain
 improvements of the general lower bounds on $\Ts(n,F)$
 given in Proposition \ref{p:turlo},
  for those $n$ which are not divisible by $q-1$.

\section{Concluding remarks and open problems}

There are several interesting directions deserving further study, which we only
 indicate briefly here.
In fact some of the preceding results can be directly extended in one way or another,
 but a more systematic study would be necessary beyond pure generalizations.

\paragraph{The quadratic bound.}

 We have seen that $(\R{G}  - 1)^2+1$ is an easy upper bound on $\Rs(G)$.
On the other hand, from the graphs studied here it seems that this naive bound is not very bad.
 In this direction we propose the following conjecture.

\bcj ~~   \label{cj:sq}
  \begin{itemize}
    \item[$(i)$] (weak form) \
 There exists a constant $c>0$ such that
  $$
    \Rs(G) \ge c \, (\R{G})^2
  $$
   holds for all graphs $G$.
    \item[$(ii)$] (strong form) \
   If $G_1,G_2,\dots$ is an infinite sequence of graphs
   without isolated vertices,
     then $\Rs(G_n) = (1-o(1)) \, (\R{G_n})^2$
     as $n\to\infty$.
  \end{itemize}
\ecj

\brm
 Proposition \ref{p:isol} implies the validity of part $(i)$
  for graphs containing very many isolated vertices.
 On the other hand the same proposition indicates that
  part $(ii)$ needs the exclusion of isolates --- or at least some
  related condition --- because otherwise $\Rs(G_n)$ is
  quadratic in $|V(G_n)|$ rather than in $\R{G_n}$.
\erm

\paragraph{More than two colors.}

Instead of 2-coloring the edges of $K_m$ one may consider $t\ge 3$ colors.
In this case the notion of singular subgraph may be introduced in several ways;
 here we mention those two of them which can be considered weakest and strongest.
In both of them we assume that $t$ graphs $G_1,\dots,G_t$ have been specified;
 moreover in any edge $t$-coloring of $K_m$ we consider the graphs $F_1,\dots,F_t$
 where the edge set of $F_j$ consists of the edges colored $j$.
Let us introduce the following notions.
\begin{itemize}
  \item A monochromatic subgraph $H$ of color $i$ is \textit{weakly singular} if
    $V(H)$ is singular in $F_i$.
  \item A monochromatic subhraph $H$ is \textit{strongly singular} if
    $V(H)$ is singular in $F_j$ for all $1\le j\le t$.
\end{itemize}
Then the weak singular Ramsey number $\Rs_w (G_1,\dots,G_t)$ is the smallest integer $n$
 such that, for every $m\ge n$, every edge $t$-coloring of $K_m$ contains a weakly
 singular subgraph $G_i$ in the color class $i$ for some $1\le i\le t$;
 and the strong singular Ramsey number $\Rs_s (G_1,\dots,G_t)$ is defined analogously.

It can be proved in various ways that $\Rs_w$ and $\Rs_s$ are finite whenever the
 graphs $G_i$ are finite.
As $t$ grows, there is an increasing number of possibilities to introduce notions
 between weak and strong singularity; and in general we have
 $\Rs_w (G_1,\dots,G_t) \le \Rs_s (G_1,\dots,G_t)$
  for all $t$ and all choices of the $G_i$.

We expect that $\Rs_w$ can be estimated more tightly than $\Rs_s$.
With the notation $n_j = (|V(G_j)|-1)^2+1$, a simple argument
 similar to the proof of Theorem \ref{th:1} yields
$$
  \Rs_w (G_1,\dots,G_t) \le \R{K_{n_1},\dots,K_{n_t}}
$$
 but this is probably quite far from being sharp in general.
For small graphs $G_i$, however, perhaps the upper bound is not terribly large.
In particular, the inequality implies $\Rs_w (K_3,K_3,K_3)\le \R{K_5,K_5,K_5}$.

\bpm
 Determine $\Rs_w (K_3,K_3,K_3)$.
\epm

 The $k$-singular generalization may also be worth studying.
For instance, in a way as in Theorem \ref{th:1},
 one can easily see that
$$
  \Rs_w (G_1,\dots,G_t, k) \le \R{K_{n_1(k)},\dots,K_{n_t(k)}}
$$
 where $n_j(k) = k(|V(G_j)|-1)^2+1$

\paragraph{Some simple graphs.}

There are some classes of graphs for which the Ramsey number is known.
 For example, one may consider

\bpm
 Determmine $\Rs(G)$ for
  \begin{itemize}
    \item[$(i)$] $G=tK_2$,
    \item[$(ii)$] $G=P_t$,
    \item[$(iii)$] $G=C_t$,
  \end{itemize}
   for all values of $t$.
\epm

\paragraph{The \boldmath{$k$}-singular version.}

So far we have a tight result concerning $k$-singular Ramsey numbers only for $P_3$
 and its subgraphs (and for edgeless graphs).
On the other hand, some estimates can easily be extended in this direction
 (cf.\ Corollary \ref{c:square}(i) and Theorem~\ref{th:2}).
It would be interesting to see the general effect of $k$ on the behavior of $\Rs(G,k)$,
 or at least for some particular examples of $G$.

\paragraph{Isolated vertices.}

We have shown in  Proposition \ref{p:isol} that the quadratic lower bound is tight
 whenever the number of non-isolated vertices is rather small compared to
  the order of the graph.
Motivated by this, the following problem is of interest.

\bpm
 Given a graph $G$, determine
 the minimum number of isolated vertices which should be added to $G$
 so that the obtained graph $G^+$ satisfies the equality $\Rs(G^+)=(|V(G^+)|-1)^2 + 1$.

\epm

\paragraph{Other structures.}

Ramsey theory has been studied for various structures, and the notion of singularity
 can be extended in a meaningful way in some of them.
For example, for any family ${\cF}$ of hypergraphs and for every natural
 number $k$, the inequality $\Rs({\cF},k) \le k(\R{\cF} - 1)^2 +1$ of
  Theorem \ref{th:1} remains valid.

\paragraph{Singular Tur\'an numbers.}

We have determined tight asymptotics for $\Ts(n,F)$ for all graphs $F$
 having at least one edge, but the exact value is known in a
 small number of cases only.
This leaves several interesting problems open.

\bpm
 Determine $\Ts(n,K_3)$ for $n\not\equiv 2 \ (\mbox{\rm mod } 4)$.
\epm

Let us note that the upper bounds in Proposition \ref{p:trian} are
 tight for $n=4$ and $n=5$, while it seems plausible to guess that
 for every other $n$ divisible by 4 the lower bound of $(i)$ gives
 the correct value.
In the other cases the lower bounds may turn out to be tight,
 at least asymptotically.

\bpm
 Determine $\Ts(n,C_4)$.
\epm

\bpm
 Prove or disprove: If $p\ge 3$ is fixed and $m$ is sufficiently large,
  then the complete $(p-1)^2$-partite graph, in which each of
  $m,m+1,\dots,m+p-2$ is the size of exactly $p-1$ vertex classes,
  is extremal for singular $K_p$, i.e.\ has $\Ts(n,K_p)$ edges
  where $n$ is the corresponding number of vertices, namely for
   $n=(p-1)\cdot\sum_{i=0}^{p-2} (m+i) = (p-1)^2\cdot (m + p/2 -1)$.
\epm

\bcj
 For every graph $F$ with $p$ vertices and chromatic number $q$,
  and every residue class $r\neq 0$ modulo $q-1$,
  there exists a constant $c(F,r)$ such that
   $$
     \lim _{n\to\infty \, ; \ n\equiv r \ (\mbox{\scriptsize\rm mod } (q-1))}
       \frac{\mbox{\rm ex}(n,K_{(p-1)(q-1)+1}) - \Ts(n,F)}{n} = c(F,r).
   $$
\ecj

\bpm
 Determine the value of the constants $c(F,r)$ for particular classes
  of graphs $F$, including complete graphs, complete bipartite graphs,
   paths and cycles.
\epm

\bpm
 Given a constant $c$ in the range $0<c<1$, find tight asymptotics on
  $\Ts(n,K_{cn})$.
\epm

\bpm
 Given $F$, find tight asymptotics on the $k$-singular Tur\'an numbers
  $\Ts(n,F,k)$.
\epm

\paragraph{Acknowledgements.}

We are grateful to the referees for their careful reading and helpful advices
 that improved the presentation of the paper.
Research of the second author was supported in part by the National Research,
 Development and Innovation Office -- NKFIH under the grant SNN 129364.

\end{document}